\RequirePackage{fix-cm}
\documentclass[smallextended]{svjour3}       
\smartqed  
\usepackage{graphicx}
\usepackage[authoryear]{natbib}
\usepackage[active]{srcltx}
\usepackage{amssymb
,amsmath,amsfonts,bbm,pifont,upgreek,bbold,
} 
\usepackage{color}
\font\teneufm=eufm10
\font\seveneufm=eufm7
\font\fiveeufm=eufm5
\newfam\eufmfam
\textfont\eufmfam=\teneufm
\scriptfont\eufmfam=\seveneufm
\scriptscriptfont\eufmfam=\fiveeufm
\newcommand\beq[1]{ \begin{equation}\label{#1} }
\newcommand{\eeq}{ \end{equation} }
\newcommand{\beqno}{ \[ }
\newcommand{\eeqno}{ \] }
\newcommand\beqa[1]{ \begin{eqnarray} \label{#1}}
\newcommand{\eeqa}{ \end{eqnarray} }
\newcommand{\beqano}{ \begin{eqnarray*} }
\newcommand{\eeqano}{ \end{eqnarray*} }
\newcommand\arr[1]{\left\{\begin{array}{l}#1\end{array}\right.}
\renewcommand{\theequation}{\arabic{section}.\arabic{equation}}
\newcommand\equ[1]{{\rm (\ref{#1})}}

\expandafter\chardef\csname pre amssym.def
at\endcsname=\the\catcode`\@
\catcode`\@=11
\def\undefine#1{\let#1\undefined}
\def\newsymbol#1#2#3#4#5{\let\next@\relax
 \ifnum#2=\@ne\let\next@\msafam@\else
 \ifnum#2=\tw@\let\next@\msbfam@\fi\fi
 \mathchardef#1="#3\next@#4#5}
\def\mathhexbox@#1#2#3{\relax
 \ifmmode\mathpalette{}{\m@th\mathchar"#1#2#3}%
 \else\leavevmode\hbox{$\m@th\mathchar"#1#2#3$}\fi}
\def\hexnumber@#1{\ifcase#1 0\or 1\or 2\or 3\or 4\or 5\or 6\or 7\or
8\or
 9\or A\or B\or C\or D\or E\or F\fi}
\ifcase\@ptsize
 \font\tenmsb=msbm10
 \font\sevenmsb=msbm7
 \font\fivemsb=msbm5
\or
 \font\tenmsb=msbm10 scaled \magstephalf
 \font\sevenmsb=msbm7 scaled \magstephalf
 \font\fivemsb=msbm5  scaled \magstephalf
\or
 \font\tenmsb=msbm10 scaled \magstep1
 \font\sevenmsb=msbm7 scaled \magstep1
 \font\fivemsb=msbm5 scaled \magstep1
\fi
\newfam\msbfam
\textfont\msbfam=\tenmsb
\scriptfont\msbfam=\sevenmsb
\scriptscriptfont\msbfam=\fivemsb
\edef\msbfam@{\hexnumber@\msbfam}
\def\Bbb#1{\fam\msbfam\relax#1}
\def\widehat#1{\setboxz@h{$\m@th#1$}%
 \ifdim\wdz@>\tw@ em\mathaccent"0\msbfam@5B{#1}%
 \else\mathaccent"0362{#1}\fi}
\def\widetilde#1{\setboxz@h{$\m@th#1$}%
 \ifdim\wdz@>\tw@ em\mathaccent"0\msbfam@5D{#1}%
 \else\mathaccent"0365{#1}\fi}

\def\RIfM@{\relax\ifmmode}
\def\nonmatherr@#1{\errmessage{\string#1\space allowed only in math mode}}
\def\Bbb{\RIfM@\expandafter\Bbb@\else
 \expandafter\nonmatherr@\expandafter\Bbb\fi}
\def\Bbb@#1{{\Bbb@@{#1}}}
\def\Bbb@@#1{\fam\msbfam\relax#1}
\def\setboxz@h{\setbox\z@\hbox}
\def\wdz@{\wd\z@}
\catcode`\@=\csname pre amssym.def at\endcsname
\newcommand{\etc}{{\it etc.\,}}

\newcommand{\nl}{{\smallskip\noindent}}
\newcommand{\noi}{{\noindent}}

\newcommand{\dst}{\displaystyle}
\newcommand\ovl[1]{ \overline {#1} }

\newcommand\su[1]{ \frac{1}{ {#1}} }

\newcommand{\torus}{ {\Bbb T}   }
\renewcommand{\natural}{ {\Bbb N}   }
\newcommand{\real}{ {\Bbb R}   }
\newcommand{\integer}{ {\Bbb Z}   }
\newcommand{\complex}{ {\Bbb C}   }
\newcommand{\rational}{ {\Bbb Q}  }

\renewcommand{\a }{ {\alpha}   }
\renewcommand{\b}{ {\beta}   }
\newcommand{\g}{ {\gamma}   }
\newcommand{\G}{ {\Gamma}   }
\renewcommand{\d}{ {\delta}   }
\newcommand{\D}{ {\Delta}   }

\renewcommand{\L}{ {\Lambda}   }
\newcommand{\m}{ {\mu}   }
\newcommand{\n}{ {\nu}   }

\newcommand{\p}{ {\pi}   }
\renewcommand{\P}{ {\Pi}   }
\renewcommand{\r}{ {\rho}   }
\newcommand{\s}{ {\sigma}   }

\renewcommand{\t}{ {\tau}   }
\newcommand{\f}{ {\varphi}   }

\renewcommand{\o}{ {\omega}   }

\newcommand{\cA}{ {\cal A} }
\newcommand{\cB}{ {\cal B} }
\newcommand{\cE}{ {\cal E} }

\newcommand{\cH}{ {\cal H} }
\newcommand{\cK}{ {\cal K} }
\newcommand{\cC}{ {\cal C} }
\newcommand{\cD}{ {\cal D} }

\newcommand{\cM}{ {\cal M} }

\newcommand{\cP}{ {\cal P} }
\newcommand{\cI}{ {\cal I} }
\newcommand{\cJ}{ {\cal J} }

\newcommand{\cS}{ {\cal S} }

\newcommand\ppu{{ (1) }}
\newcommand\ppd{{ (2) }}

\newcommand\pph{{ (h) }}

\newcommand\ppj{{ (j) }}

\newcommand\ppi{{ (i) }}

\newcommand\ppo{{ (0) }}

\newcommand\GG{{\rm G}}

\newcommand\id{{\, \rm id \,}}

\begin{document}

\title{A first integral to the partially averaged Newtonian potential of the  three--body problem
\thanks{
I wish to thank the Associate Editor and the anonymous Reviewers, who, with their constructive remarks, considerably helped me to improve the presentation of the results. I am indebted to
 Giacomo Tommei, for providing me  the reference~\citet{bekovO78}, which was the starting point of this study. I heartily thank  A. Celletti, C. Efthimiopulos, Marcel Guardia, Tere Seara and Lei Zhao for their interest.
I am grateful to the Mathematical Sciences Research Institute in Berkeley for its kind hospitality in Fall 2018, during which the paper was written.\\
 This research has been funded by the European Research Council [Grant 677793 Stable and Chaotic Motions in the Planetary Problem].
 }
}
%



\author{Gabriella Pinzari        
}


\institute{G. Pinzari \at
              Dipartimento di Matematica ``T. Levi--Civita'' \\
              via Trieste, 63, 
              35131, Padova\\
              Tel.: +39-049-8271477\\
              \email{gabriella.pinzari@math.unipd.it}    \\
              ORCID 0000-0002-0031-7046
}

\date{Received: / Accepted: date}

\maketitle

\begin{abstract}
We consider the partial average i.e., the Lagrange average with respect to {\it just one} of the two mean anomalies,
of the  Newtonian part of the perturbing function in the three--body problem Hamiltonian. We prove that such a partial average exhibits a   non--trivial first integral. We  show that this integral is fully responsible of certain cancellations in the averaged Newtonian potential, including a property noticed by Harrington in the 60s. We also highlight its joint r\^ole (together with certain symmetries) in the appearance of the so called ``Herman resonance''. Finally, we discuss an application and an open problem. 
\end{abstract}

\section{Motivation}
\setcounter{equation}{0}
\renewcommand{\theequation}{\arabic{equation}}

The  purpose of this work is to highlight a  property of the ``partial average of the Newtonian potential'' and discuss some consequences.

\nl
{By ``partial averaged Newtonian potential'' we mean the following. Let $(y^\ppi, x^\ppi)$ $=$ $\big((y^\ppi_1$, $y^\ppi_2$, $y^\ppi_3)$, $( x^\ppi_1$, $x^\ppi_2$, $x^\ppi_3)\big)$, with $i=1$, $2$, be impulse--position coordinates for a  two--particle system (which we also call ``planets'') and let
\beqa{C}\cC:\qquad(\L_2,\ell_2,u,v)\in \cA\times\torus\times V\to (y, x)=(y^\ppu, y^\ppd, x^\ppu, x^\ppd)\in (\real^3)^4\eeqa
 where $\cA$ is a domain\footnote{{By ``domain'' we mean an open and connected set in ${\mathbb K}=\real^m, \complex^m$.}} in $\real$,   $V$ is a domain in $\real^{10}$, $\torus:=\real/(2\p\integer)$, $(u, v)=\big((u_1, u_2, u_3, u_4, u_5)$, $(v_1, v_2, v_3, v_4, v_5)\big)$, 
be a change of coordinates, which we call, for brevity, {\it partial Kepler map},  which ``preserves the standard two--form'':
$$dy^\ppu\wedge dx^\ppu+dy^\ppd\wedge dx^\ppd=d\L_2\wedge d\ell_2+d u\wedge d v$$
and ``integrates the Keplerian motions of $(y^\ppd, x^\ppd)$'':
 \beq{2B}
\left(\frac{\|y^\ppd\|^2}{2{\rm m}_2}-\frac{{\rm m}_2{{\rm M}}_2}{\|x^\ppd\|}\right)\circ\cC=
-\frac{{\rm m}_2^3{\rm M}_2^2}{2\L_2^2}=:{\rm h}_{\rm Kep}^\ppd(\L_2)\ ,\eeq 
where ${\rm m}_2$, ${\rm M}_2$ are suitable ``mass parameters''. Of course, we have assumed that the image of $\cC$ in~\equ{C}
is a domain of  $(y, x)$ where the left hand side of~\equ{2B} takes negative values.  We also assume, throughout the paper, that $(y, x)$ are chosen so that the instantaneous ellipse ${\mathbb E}_2$ generated by the two--body Hamiltonian~\equ{2B} has non--vanishing\footnote{For simplicity, we refrain to  formulate the results in the case that the map $\cC$ in~\equ{C} is regular when the eccentricity of ${\mathbb E}_2$ vanishes, as it happens, for example, in the case of the Poincar\'e or the {\sc rps} map.} eccentricity, so we
denote as ${\rm P}^\ppd$, $\|{\rm P}^\ppd\|=1$, the direction of its perihelion. The  angle $\ell_2$ will be referred to as ``mean anomaly'', for uniformity with the name attributed to an analogue angle in the set of the coordinates named after Delaunay (see   e.g.~\citet{fejoz13} for a definition).
 We look at the Lagrange average
 \beqa{h1} h_2(\L_2,u,v):=\frac{1}{2\p}\int_\torus \frac{d\ell_2}{\|x^\ppu(\L_2,\ell_2, u, v)-x^\ppd(\L_2,\ell_2, u, v)\|}\eeqa which we will refer to as  {\it partially averaged Newtonian potential}.
\vskip.1in
\noi
There are many examples, in Celestial Mechanics, of canonical maps of the form above. Well known ones are the above mentioned Delaunay map (hereafter, $\cD$), or the  coordinates after the Jacobi--Deprit reduction\footnote{{The coordinates discovered by~\citet{deprit83} are an extension,  to any number of particles, of~\citet{jacobi1842}, which hold only  for a two--particle system.}} of the nodes ($\cJ$)~\citep{jacobi1842, deprit83}. Another example, called ``perihelia reduction'' ($\cP$), has been introduced by~\citet{pinzari18}.  A comprehensive  review can be found in~\citet{pinzari14}.
All the maps mentioned here might actually be named  {\it double Kepler maps}, since, in such cases, they satisfy~\equ{C}--\equ{2B}, with, in turn, the  $(u, v)$'s having the form
$$(u, v)=(\L_1, \ell_1, \hat u, \hat v)\ ,\quad du\wedge dv=d{\L_1}\wedge d\ell_1+d\hat u\wedge d\hat v$$
where $\ell_1\in \torus$ and  $\L_1$ is such that~\equ{2B} holds also with ${\rm m}_2$, ${\rm M}_2$, $y^\ppd$, $x^\ppd$ replaced by ${\rm m}_1$, ${\rm M}_1$, $y^\ppu$, $x^\ppu$.  In Section~\ref{The partial Kepler map} below we shall present a ``genuine'' partial Kepler map, namely a map $\cC$ where~\equ{2B} holds only for one of the bodies. \vskip.1in 
\noi
We have been interested to  the function~\equ{h1} because, in  planetary $(1+N)$-body theories, 
one has to deal with analogue
maps of  the kind
$$\cC_N:\quad (\L, \ell, \hat u, \hat v)\in \cA^N\times \torus^N\times W\to \big((y^\ppu, \cdots, y^{(N)}), (x^\ppu, \cdots, x^{(N)})\big)$$
with $W$ a domain in $\real^{4N}$ in terms of which the Hamiltonian of the system is
\beqa{planetary} H_N(\L, \ell,  \hat u, \hat v)=-\sum_{i=1}^N \frac{{\rm m}_i^3{\rm M}_i^2}{2\L_i^2}+f_N(\L, \ell, \hat u, \hat v)\eeqa
where the Hamiltonian is composed of a leading``Keplerian part'', given by $-\sum_{i=1}^N \frac{{\rm m}_i^3{\rm M}_i^2}{2\L_i^2}$, 
slightly perturbed by a  function  $f_N$.
{The  splitting \equ{planetary} is possible -- and in fact it has been adopted by~\citet{arnold63, laskarR95, chierchiaPi11b,  palacianSY12, palacian18, pinzari18}, for example, in the so--called {\it planetary problem}  where one of the masses (``sun'') is much larger than the remaining, equally  sized, $N$ ones (``planets'').  {In that case,} averaging over the Keplerian frequency vector \beqa{kep freq}\o_{\rm Kep}=(\o_{\rm Kep, 1}, \cdots, \o_{{\rm Kep}, N})\ ,\quad \o_{\rm Kep, i}=\frac{{\rm m}_i^3{\rm M}_i^2}{\L_i^3}\eeqa   leads to study the so--called ``secular problem''
$$\ovl H_N(\L, \hat u, \hat v)=-\sum_{i=1}^N \frac{{\rm m}_i^3{\rm M}_i^2}{2\L_i^2}+\ovl f_N(\L, \hat u, \hat v)$$
 where the  perturbing term  is given by  ``multi--averaged Newtonian potential'' (the study of which goes back to ~\citet{sundman1916})
 \beqa{multiav} \ovl f_N(\L,\hat u, \hat v)&:=&{-}\frac{1}{(2\p)^N}\sum_{1\le i<j\le N}m_i m_j\nonumber\\
 &&\int_{\torus^N} \frac{d\ell_id\ell_j}{\|x^\ppi(\L_i,\L_j,\ell_i, \ell_j, \hat u, \hat v)-x^{{\ppj}}(\L_i,\L_j,\ell_i, \ell_j, \hat u, \hat v)\|}\ .\eeqa
It is known~\citep{gallavotti86} that  the dynamics of the full problem is well approximated by the one of the secular one as soon as no resonances between the frequencies~\equ{kep freq} appear. In case of resonance, for example, in the case $N=2$, with  the two planets being much distant one to the other, it is reasonable to expect that a better approximation is obtained replacing the  average~\equ{multiav} with  the partial average~\equ{h1}. }
{Concretely,  it might be  challenging to investigate whether there is an application  to any of the following regions of motion that  have been proposed   by~\citet[p. 310]{fejoz02}, for the $N=2$ case:
\begin{itemize}
\item[--] the {\it planetary} region, where the eccentricity of the outer ellipse and both semi--major axes are in a small compact set, and two masses are small compared to the third mass;
\item[--] the {\it lunar}  region, where the masses are in a compact set, and the outer body is far away from the outer two;
\item[--] the {\it anti--planetary} region, where the outer body ellipse may have  a large mass, provided its ellipse is far away from the outer two;
\item[--] the {\it anti--lunar} region, when the  ellipses  of the two outer bodies are close, but the corresponding masses are much different.
\end{itemize}
}

\nl
We now go back to $h_2$ in~\equ{h1}. We firstly observe that

\begin{theorem}\label{theo: Harr} $h_2$ is integrable by quadratures.
\end{theorem}

\nl
Indeed, $h_2$ has  six degrees of freedom and possesses, besides itself,  the following five commuting\footnote{In Hamiltonian mechanics $f(p, q)$, $g(p, q)$ are said to be {\it Poisson--commuting} if their {\it Poisson parentheses} $\{f, g\}:=\sum\partial_p f\partial_q g-\partial_p g\partial_q f$ vanish. Poisson commutation of $f$ and $g$ is equivalent to say that $g$ remains constant along the Hamiltonian motions of $f$. } integrals: 

\nl
 $\rm I_1$:= the semi--major axis action $\L_2:={\rm m}_2\sqrt{{\rm M}_2 a_2}$;\\
$\rm I_2$:= the Euclidean length $\|x^\ppu\|$ of $x^\ppu$;\\
$\rm I_3$:=   the Euclidean length of  the  total angular momentum ${\rm C}:={\rm C}^\ppu+{\rm C}^\ppd$, with ${\rm C}^\ppi:=x^\ppi\times y^\ppi$, and ``$\times$'' denoting skew--product;\\
$\rm I_4$:= its third component; \\
$\rm I_5$:= the projection of the angular momentum ${\rm C}^\ppd$ along the direction $x^\ppu$.

\nl
{Indeed, $\rm I_1$ is trivially due to the $\ell_2$--averaging; $\rm I_3$ and $\rm I_4$ descend from the invariance by rotations of $h_2$; $\rm I_2$ and $\rm I_5$ from invariance by rotations around the $x^\ppu$ axis.}
{Such integrals are independent if ${\rm C}^\ppu$ and ${\rm C}^\ppd$ are not parallel. Otherwise, the problem reduces to be planar, namely, $h_2$ has four degrees of freedom, and three independent commuting integrals are obtained neglecting, in the list above,  ${\rm I}_4$ and ${\rm I}_5$.}

\begin{remark}\rm 
The list of independent first integrals to $h_2$ is even longer then the one above. For example, in the spatial case, the three components of  $x^\ppu$ and the three components ${\rm C}^\ppd$ are {\it all} first integrals. However, the maximum number of {\it commuting} first integrals that can be formed with these quantities is four (and the functions $\rm I_2$, $\rm I_3$, $\rm I_4$ and $\rm I_5$ are an example of them).
\end{remark}

\begin{remark}\rm 
{The integrability of $h_2$ does not imply that also the partial average of the three--body problem  Hamiltonian  is so, because this includes also a kinetic term. This is an even different situation compared to the secular problem mentioned above, whose non--integrability  is clearly proven, as a consequence of the so--called {\it splitting of separatrices}~\citep{fejozG15}.}
\end{remark}


\nl
{We now consider the ellipse generated by the ``Kepler Hamiltonian'' at left hand side in~\equ{2B} and denote as  ${\rm e}_2:=\sqrt{1-\frac{{\rm G}_2^2}{\L_2^2}}$ its eccentricity, where ${\rm G}_2:=\|{\rm C}^\ppd\|$. Then, let
\beqa{cG}{{\rm E}_0}:={\rm G}_2^2-{\rm m}_2^2{\rm M}_2 {\rm e}_2\,x^\ppu\cdot {\rm P}^\ppd\eeqa
The following fact is a bit more subtle. }

\begin{theorem}\label{partial integral}
The function ${\rm E}_0$ is a first integral of $h_2$. 
\end{theorem}
\proof The proof of this theorem uses some results from\footnote{${\rm m}_2$, ${\rm M}_2$, ${\rm M}_1$,  correspond to ${\rm m}$, $\cM$, $\m\cM$ in~\citet{pinzari18a}.}~\citet{pinzari18a}, that here we recall. We consider the Hamiltonian
$${\rm J}=\frac{\|y^\ppd\|^2}{2{\rm m}_2}-\frac{\rm m_2{\rm M}_2}{\|x^\ppd\|}-\frac{{\rm m}_2{\rm M}_1}{\|x^\ppu-x^\ppd\|}\ .$$
This is the Hamiltonian of one moving particle $(y^\ppd, x^\ppd)$ having mass ${\rm m}_2$, subject to the gravitational attraction by two fixed particles: ${\rm M}_2$, at the origin, and ${\rm M}_1$, at $x^\ppu$. The Hamiltonian is integrable by quadratures, for having, as first integrals, the function  {${\rm I}_5$} defined above (which trivialises in the case of the planar problem) and the function
$${\rm E}={\rm E}_0+{\rm M}_1 {\rm E}_1$$
where ${\rm E}_0$ is as in~\equ{cG}, while
$${\rm E}_1={\rm m}_2^2\frac{x^\ppu\cdot(x^\ppu-x^\ppd)}{\|x^\ppu-x^\ppd\|}\ .$$
We write ${\rm J}$ and ${\rm E}$ in terms of a given partial Kepler map, $\cC$. We obtain
\beqa{Jc}{\rm J}_\cC=-\frac{{\rm m}_2^3{\rm M}_2^2}{2\L_2^2}-\frac{{\rm m}_2{\rm M}_1}{\|x_\cC^\ppu-x_\cC^\ppd\|}\ ,\qquad {\rm E}_\cC={\rm E}_{0,\cC}+{\rm M}_1{\rm E}_{1, \cC}\eeqa
where the subfix $\cC$ denotes the composition with $\cC$. The commutation  of ${\rm J}_\cC$ and  ${\rm E}_\cC$ implies the following relation, which is obtained picking up the terms at the first order in ${\rm M}_1$:
$$\left\{-\frac{{\rm m}_2^3{\rm M}_2^2}{2\L_2^2},\ {\rm E}_{1, \cC}\right\}+\left\{-\frac{{\rm m}_2}{\|x_\cC^\ppu-x_\cC^\ppd\|},\ {\rm E}_{0, \cC}\right\}=0	 \ .$$
Taking the $\ell_2$--average of this identity, the first term vanishes itself:
$$\frac{1}{2\p}\int_\torus \left\{-\frac{{\rm m}_2^3{\rm M}_2^2}{2\L_2^2},\ {\rm E}_{1, \cC}\right\} d\ell_2=\frac{1}{2\p}\frac{{\rm m}_2^3{\rm M}_2^2}{\L_2^3}\int_{\torus}\partial_{\ell_2}{\rm E}_{1, \cC}d\ell_2\equiv0\ .$$
Hence,
$$0=\frac{1}{2\p}\int_\torus \left\{-\frac{{\rm m}_2}{\|x_\cC^\ppu-x_\cC^\ppd\|},\ {\rm E}_{0, \cC}\right\}d\ell_2=
\big\{-{\rm m}_2\, h_2,\ {\rm E}_{0, \cC}\big\}$$
since ${\rm E}_{0, \cC}$ is $\ell_2$--independent. This is the thesis. \qed

\vskip.1in
\noi
In the next sections  we highlight some properties of the partially averaged Newtonian potential that descend from  Theorem~\ref{theo: Harr} and Theorem~\ref{partial integral}. More precisely,  the paper is organised as  follows. In Section~\ref{Generalised Harrington property} we show that, as a consequence of Theorem~\ref{partial integral}, an infinite number of Fourier coefficients in the expansion of $h_2$ with respect to the perihelion of its outer planet cancel. This property is a generalisation of a fact noticed by S. Harrington in the 1960s~\citet{harrington69}. To this purpose, we introduce a set of canonical coordinates in terms of which $h_2$ and ${\rm E}_0$ are reduced to one degree of freedom. In Section~\ref{Renormalizable integrability} we show that there is an explicit functional dependence between $h_2$ and ${\rm E}_0$. We call this circumstance  ``renormalizable integrability''. 
{The author argues that it might be helpful in the framework of the study of the three--body problem. For example, 
it would be nice to understand if
fixed points of ${\rm E}_0$, both of elliptic and hyperbolic character, being at the same time fixed points to $h_2$ with the same character, might give rise to quasi--periodic motions in the three--body problem; if hyperbolic equilibria might lead to a splitting of separatrices, \etc. Instead of addressing such  issues here (which would lead much further than the purposes of this note; see, however,~\citet{pinzari18a} for an application in this direction),  we discuss the relations between level curves and the fixed points of the two functions.} 
Next, we show that, as a consequence of  {renormalizable integrability} and the well known 
\begin{proposition}[Keplerian property]\label{prop: kepler}
\beqano\frac{1}{2\p}\int_{\torus}\frac{d\ell_2}{\|x_\cC^\ppd\|}=\frac{1}{a_2}\quad \forall\ \cC\ .\eeqano
\end{proposition}
A linear combination with integer coefficients in a suitable expansion of $h_2$ is identically verified. We  name it ``generalised Herman resonance'' since it recalls the well known Herman resonance in the doubly averaged Newtonian potential (we refer to~\citet{abdullahA01} or~\citet[Propriet\'e 80]{fejoz04} for informations on Herman resonance).
{After proving, in} Section~\ref{A algebraic property of Legendre polynomials}, an  algebraic property of the well known {Legendre} polynomials (which, roughly, says that a certain average of a {Legendre} polynomial is still a {Legendre} polynomial),  we establish, {in Section~\ref{Applications}}, a link between the aforementioned generalised Herman resonance and Herman resonance.  In this conclusive section we also provide a sort of ``eccentricity--inclination'' expansion at any order for such function and discuss a problem which is left open.

\section{Generalised Harrington property}\label{Generalised Harrington property}

In this section we assume that the map $\cC$ in~\equ{C}  includes, among the $u$'s, the impulse $u_1:={\rm G}_2:=\|{\rm C}^\ppd\|$. We also give  $x^\ppu$, $x^\ppd$ the meaning of ``interior'', ``exterior'' planet, respectively,  because we write formal expansions with respect to $\|x^\ppu\|$.

\nl
We prove the following
\begin{theorem} \label{HP} 
Fix a domain for $\cC$ where $x_\cC^\ppu\times {\rm C}_\cC^\ppd$, ${\rm C}_\cC^\ppd\times{\rm P}^\ppd_\cC$, and ${\rm C}_\cC^\ppd$ never vanish.
Let 
 $$h:\quad (\L_2, u, v)\in \cA\times V\to h(\L_2, u, v)$$
Poisson--commute with ${\rm E}_0$. Assume that
$h$ has the form
\beqa{legendre exp}h=\sum_{n=0}^{\infty}\sum_{m=0}^{+\infty} h_{nm}(\L_2, u, v)\r^n \cos{m\f}\ ,\eeqa
where $\r(\L_2, u, v):=\|x^\ppu_\cC\|$  and $\f (\L_2, u, v)$ is the angle formed by the two vectors $x_\cC^\ppu\times {\rm C}_\cC^\ppd$, ${\rm C}_\cC^\ppd\times{\rm P}^\ppd_\cC$, with respect to the counterclockwise  orientation established by ${\rm C}_\cC^\ppd$. Assume also that $h_{nm}$ depends on $(\L_2, u, v)$ only via the following quantities
\beqa{quant***}\L_2, \quad u_1={\rm G}_2\ ,\quad \Theta:=\frac{x^\ppu_\cC\cdot {\rm C}^\ppd_\cC}{\|x^\ppu_\cC\|}\ ,\eeqa
with
 $h_{0m}$ being independent of $u_1={\rm G}_2$ for all $m\ge 0$. Then
 \beqa{Hprop} h_{nm}(\L_2, u, v)\equiv 0\qquad {\rm if}\qquad m\ge \max\{1,\ n\}\ ,\ \forall\ n\ge 0\ .\eeqa
 In the case that $h_{nm}=0$ for $n-m$ odd, for $n\ge1$, the following stronger identities hold:
 \beqa{Hprop1} h_{nm}(\L_2, u, v)\equiv 0\qquad {\rm if}\qquad m\ge n-1\ ,\quad \forall \ n\ge 1\ .\eeqa
  \end{theorem}
To prove Theorem~\ref{HP}, we shall need the following
\begin{lemma}\label{lem: precise rules}
Let the  functions
$$h(\G,\g)=\sum_{n=0}^\infty\sum_{m=0}^\infty\varepsilon^n h_{nm}(\G)\cos m\g\qquad g(\G,\g)=a(\G)+\varepsilon b(\G)\cos\g$$
verify
\beqa{comm**}\Big\{h, g\Big\}_{\G,\g}:=\partial_\G h\partial_\g g-\partial_\G g\partial_\g h\equiv 0\eeqa
and assume that {$\partial_\G a\not\equiv 0$} and $h_{0m}$ is independent of $\G$ for all $m\ge 0$. 
Then $h_{nm}=0$ for all $m\ge \max\{1, n\}$.
\end{lemma}
\proof Due to the assumptions of  $h$ and $g$, their Poisson parenthesis at left hand side of~\equ{comm**} is a Fourier series including only {sines} $\{\sin k\g\}_{k\ge 1}$.
Projecting ~\equ{comm**} over such  basis, we obtain the following relations:
\beqa{nm}
 m \partial_\G ah_{nm}&=&-\su 2\Big((m-1)h_{n-1,m-1}+(m+1)h_{n-1,m+1}\Big)\partial_\G b\nonumber\\
&+&\su 2\Big(\partial_\G h_{n-1,m-1}- \partial_\G h_{n-1,m+1}+\partial_\G h_{n-1,0}\d_{m,1}\Big)b
\eeqa
for all $n= 0$, $1$, $m=1$, $2$, $\cdots$, where $\d_{ij}$ is the Kronecker symbol, and $h_{-1,k}:=0$ for all $k\in \integer$.
We 
now prove that such relations imply $h_{nm}\equiv0$ for $m\ge \max\{1, n\}$. We proceed by steps.
\paragraph{\rm (i)} We prove $h_{0m}=0$ for $m=1,\ 2,\ \cdots$. We use~\equ{nm} with $n=0$ and $m=1,\ 2,\ \cdots$: 
\beqano
m\partial_\G a h_{0m}&=&-\su2 \big((m-1)h_{-1, m-1}+(m+1)h_{-1, m+1}\big)\nonumber\\
&+&\su2\big(\partial_\G h_{-1, m-1}-\partial_\G h_{-1, m+1}+\partial_\G h_{-1, 0}\d_{m, 1}\big)b\nonumber\\
&\equiv&0\quad m=1,\ 2,\ \cdots
\eeqano
since  $h_{-1, k}=0$ for all $k\in \integer$, as $\partial_\G a\not\equiv0$. 
\paragraph{\rm (ii)} We prove $h_{1m}=0$ for $m\ge 1$.

\nl
\quad{(ii)-a} We prove $h_{11}=0$. We use~\equ{nm} with
$n=m=1$. We obtain
\beqano
  \partial_\G ah_{11}&=&-\su 2\Big(2h_{0,2}\Big)\partial_\G b\nonumber\\
&+&\su 2\Big(\partial_\G h_{00}- \partial_\G h_{02}+\partial_\G h_{00}\d_{11}\Big)b\nonumber\\
&\equiv&0
\eeqano
since $h_{02}=0$ by (i) and $\partial_\G h_{00}=0$ by assumption. 

\nl
\quad{(ii)-b} We prove $h_{1m}=0$ for $m\ge 2$. We use~\equ{nm} with $n=1$, $m\ge 2$:
\beqano
 m \partial_\G ah_{1m}&=&-\su 2\Big((m-1)h_{0,m-1}+(m+1)h_{0,m+1}\Big)\partial_\G b\nonumber\\
&+&\su 2\Big(\partial_\G h_{0,m-1}- \partial_\G h_{0,m+1}+\partial_\G h_{0,0}\d_{m,1}\Big)b
\nonumber\\
&\equiv&0
\eeqano
because the first lines vanishes by (i), while the second vanishes because, by assumption, $\partial_\G h_{0,p}$ for all $p\ge 1$.
\paragraph{\rm (iii)}   We prove $h_{nm}=0$ for $n\ge1$ and $m\ge n$. We proceed by induction on $n$. The case $n=1$ has been done in (ii). We assume that it is true for $n\ge 1$ and prove it for $n+1$. We use~\equ{nm} replacing $n$ with $n+1$ and taking $m\ge n+1$:
\beqano
 m \partial_\G ah_{n+1, m}&=&-\su 2\Big((m-1)h_{n,m-1}+(m+1)h_{n,m+1}\Big)\partial_\G b\nonumber\\
&+&\su 2\Big(\partial_\G h_{n,m-1}- \partial_\G h_{n,m+1}+\partial_\G h_{n,0}\d_{m,1}\Big)b\nonumber\\
&\equiv&0
\eeqano
Here we have used that for $m\ge n+1$, $m+1>m-1\ge n$, so the first line and the two first terms in the second line vanish. The last term also vanishes because $m\ge n+1\ge 2$, so the Kronecker symbol is zero. The lemma is completely proved. \qed
\vskip.1in \noi We now proceed to prove Theorem~\ref{HP}. To this end, we introduce a specific  system of canonical coordinates 
which will allow us to apply the lemma above.
\paragraph{The  $\cK$--map}\label{The partial Kepler map} 
Define the ``nodes''
$$\n_0:={\rm k}\times {\rm C}\ ,\quad \n_1:={\rm C}\times x^\ppu\ ,\quad \n_2:=x^\ppu\times {\rm C}^\ppd\ ,\quad \n_3:={\rm C}^\ppd\times {\rm P}^\ppd
$$
and assume that they do not vanish. Denote, as above,
as  ${\rm P}^\ppd$, with $\|{\rm P}^\ppd\|=1$ the direction of its perihelion (well defined because the eccentricity does not vanish), $a_2$ its semi--major axis, we  define the map
\beqno\cK:\quad \big(
\L_2,{\rm l}_{2},   {\rm Z}, {\rm G}, {\rm R}_1, {\rm G}_2, \Theta, {\rm z}, {\rm g}_2, {\rm g}, {\rm r}_1, \vartheta
\big)\to (y^\ppu_{\cK}, y^\ppd_{\cK}, x^\ppu_{\cK}, x^\ppd_{\cK})\ .\eeqno
via the relations
\beqa{K} 
\cK^{-1}:\qquad \arr{
 \dst {\rm Z}:={\rm C}\cdot {\rm k}\\ \\
 \dst {\rm G}:=\|{\rm C}\|\\ \\
  \dst  {\rm R}_1:=\frac{y^\ppu\cdot x^\ppu}{\|x^\ppu\|}\\\\
\dst  \L_2={\rm m}_2\sqrt{{\rm M}_2 a_2}\\\\
 \dst{\rm G}_2:=\|{\rm C}^\ppd\|\\\\
 \dst \Theta:=\frac{{\rm C}^\ppd\cdot x^\ppu}{\|x^\ppu\|}\\
 }\qquad\qquad \arr{
 \dst {\rm z}:=\a_{{\rm k}}({\rm i}, \n_0)\\ \\
 \dst{\rm g}:=\a_{{\rm C}}(\n_0, \n_1)\\ \\
  \dst {\rm r}_1:=\|x^\ppu\|\\ \\
{\rm l}_2:=\textrm{\rm mean anomaly of } x^\ppd\ {\rm on\ \mathbb E}\\ \\
{\rm g}_2:=\a_{{\rm C}^\ppd}(\n_2, \n_3)\\\\
  \dst\vartheta:=\a_{x^\ppu}(\n_1, \n_2)
 }
\eeqa
where $({\rm i}, {\rm j}, {\rm k})$ is a prefixed reference frame, and for $u,v\in\real^3$ lying in the plane orthogonal to a vector $w$ and $\a_w(u,v)$ denotes the positively oriented angle (mod $2\p$) between $u$ and $v$  (orientation follows  the ``right hand rule''). We remark that the planar case correspond to {taking} $\Theta=0$ and $\vartheta=\p$ (prograde case) of $\vartheta=0$ (retrograde case).
\\
The map $\cK$ verifies~\equ{C}--\equ{2B} with $\ell_2={\rm l}_2$ and $u=({\rm G}_2, \breve u)$, $v=({\rm g}_2, \breve v)$, where $\breve u=({\rm Z}, {\rm G}, {\rm R}_1, \Theta)$, $\breve v=({\rm z}, {\rm g}, {\rm r}_1, \vartheta)$. Therefore $u$ and $v$ are also as claimed in the assumptions of Theorem~\ref{HP}.
The canonical character of the coordinates $\cK$ is discussed in~\citet{pinzari18a} and to such paper we refer also for the formula,  in terms of $\cK$, of the function ${{\rm E}_0}$ in~\equ{cG}, which is
\beqa{cG***}{{{\rm E}_0}}={\rm G}_2^2+{\rm m}_2^2{\rm M}_2{\rm r}_1\sqrt{1-\frac{{\rm G}_2^2}{\L_2^2}}\sqrt{1-\frac{\Theta^2}{{\rm G}_2^2}}\cos{\rm g}_2\ .\eeqa
We continue denoting as $h$ the function in the statement expressed in terms of $\cK$. 
It follows from the definitions~\equ{K} that $\r={\rm r}_1$ and $\f={\rm g}_2$, so, by~\equ{legendre exp},  $h$ is given by 
\beqano h=\sum_{n=0}^{\infty}\sum_{m=0}^{+\infty} {\rm r}_1^n h_{nm}(\L_2, \Theta, {\rm G}_2) \cos{m {\rm g}_2}\ .\eeqano
Here we have used that, by assumption, the coefficients $h_{nm}$ in this expansion depend only on $\L_2$, ${\rm G}_2$, $\Theta$. Therefore, in terms of $\cK$, the assumption that $h$  Poisson commutes with ${{\rm E}_0}$ reduces to 
$$\big\{h, {{\rm E}_0}\big\}_{({\rm G}_2, {\rm g}_2)}=\partial_{{\rm G}_2}h\partial_{{\rm g}_2}{{\rm E}_0}-\partial_{{\rm g}_2}h\partial_{{\rm G}_2}{{\rm E}_0}\equiv0\ .$$
Furthermore, with $a={\rm G}_2^2$, we have $\partial_{{\rm G}_2} a\not\equiv 0$ and, finally,  $h_{0m}$ is independent of ${\rm G}_2$ for all $m\ge0$, being this one of the assumptions of Theorem~\ref{HP}.
We can thus apply Lemma~\ref{lem: precise rules} and we obtain that $h_{nm}({\rm r}_1, \L_2, \Theta)\equiv0$ for $m\ge \max\{1, n\}$. The identities~\equ{Hprop1} trivially follow, under the additional assumption that $h_{nm}=0$ if $n-m$ is odd. \qed
\paragraph{Application of Theorem~\ref{HP} to the function $h_2$} 
In this section we discuss the application of Theorem~\ref{HP} to the function $h_2$ in~\equ{h1}.  First of all, $h_2$ Poisson commutes with ${\rm E}_0$, as stated by Theorem~\ref{partial integral}. As in the proof of Theorem~\ref{HP}, we now write $h_2$ in terms of the coordinates
$\cK$ in~\equ{K} and we fix a domain as in the statement of the theorem. This map is  useful because $\r={\rm r}_2$, $\f={\rm g}_2$  and the functions in~\equ{quant***}  are coordinates in such system, so we have only to check that $h_2$ affords an expansion of the form:
\beqa{h1***new} h_2=\sum_{ n=0}^{\infty}\sum_{m=0}^{+\infty} {\rm r}_1^n h_{2, nm}(\L_2, \Theta, {\rm G}_2) \cos{m {\rm g}_2}\ ,\eeqa
with $ \partial_{\rm G_2}h_{1, 0m}(\L_2, \Theta, {\rm G}_2)\equiv 0$. We shall also check that, in this summand, only terms with even $n-m$ appear.
We observe that, since  $h$ commutes with ${\rm I}_1$, $\cdots$,  ${\rm I}_5$, 
and, by their definitions, such functions are coordinates in the system $\cK$:
\beqa{nice integrals}{{\rm I}_1=\L_2\ ,\quad {\rm I}_2={\rm r}_1\ ,\quad {\rm I}_3={\rm G}\ ,\quad {\rm I}_4={\rm Z}\ ,\quad {\rm I}_5=\Theta}\eeqa
then we have that $h_2$ is independent of their {conjugate} coordinates, respectively,
{$\ell_2$, ${\rm R}_1$,  ${\rm g}$, ${\rm z}$,  $\vartheta$}. The angles ${\rm g}$, ${\rm z}$ are themselves first integrals to $h_2$ and so we have that $h_2$ is also independent of ${\rm G}$, ${\rm Z}$. In summary, $h_2$ will be a function of ${\rm r}_1$, $\L_2$, $\Theta$, ${\rm G}_2$, ${\rm g}_2$ only. 
Now we check that $h_2$ affords an expansion of the form~\equ{legendre exp}, with $h_{nm}$ depending only on the quantities~\equ{quant***}. As already observed in the proof of Theorem~\ref{HP},  in terms of the coordinates $\cK$, this reduces to check that $h_2$, in terms of $\cK$ has an expansion of the form~\equ{h1***new}. To this end, we start  from 
the expansion of the Newtonian potential in Legendre  polynomials (see Section~\ref{A algebraic property of Legendre polynomials})
\beqa{Legendre}\frac{1}{\|x^\ppu-x^\ppd\|}=\sum_{n=0}^\infty \cP_n(t) \frac{\|x^\ppu\|^n}{\|x^\ppd\|^{n+1}}\qquad t:=\frac{x^\ppu\cdot x^\ppd}{\|x^\ppu\| \|x^\ppd\|}\eeqa
In terms of $\cK$, such quantities are given by
$$\|x^\ppu\|={\rm r}_1\ ,\quad{\|x^\ppd\|=a_2\frac{\frac{\GG^2_2}{\L^2_2}}{1+\sqrt{1-\frac{{\rm G}_2^2}{\L_2^2}}\cos {\rm f}_2}}\ ,\quad t=-\sqrt{1-\frac{\Theta^2}{{\rm G}_2^2}}\cos({\rm g}_2+{\rm f}_2)$$
where $a_2=\frac{\L_2^2}{{\rm m}_2^2{\rm M}_2}$;  ${\rm f}_2={\rm f}_2(\L_2, {\rm G}_2, {\rm l}_2)$ is the true anomaly. {The two former expressions are classical; the one  for $t$ has been worked out by~\citet{pinzari18a}.} Inserting these expressions into~\equ{Legendre} and taking the ${\rm l}_2$--average\footnote{{Recall the well known transition formula (see, e.g.,~\citet{palacianVY17})
$\dst d{\rm l}_2=\frac{\frac{{\rm G}_2^3}{\L_2^3}}{\left(1+\sqrt{1-\frac{{\rm G}_2^2}{\L_2^2}}\cos {\rm f}_2\right)^2}d{\rm f}_2$.
}}, we have that 
\beqa{Pn}h_2({\rm r}_1, \L_2, \Theta, {\rm G}_2, {\rm g}_2)=\sum_{n=0}^\infty h_{2, n}(\L_2, \Theta, {\rm G}_2, {\rm g}_2) {\rm r}_1^n\ . \eeqa
with
\beqa{Qn}h_{2, n}(\L_2, \Theta, {\rm G}_2, {\rm g}_2)&=&\frac{1}{2\pi a_2^{n+1}}\frac{\L_2^{2n-1}}{\GG_2^{2n-1}}\int_{\torus} 
\left(1+\sqrt{1-\frac{{\rm G}_2^2}{\L_2^2}}\cos {\rm f}_2\right)^{n-1}\nonumber\\
&&
{\cP_n\left(-\sqrt{1-\frac{\Theta^2}{{\rm G}_2^2}}\cos({\rm g}_2+{\rm f}_2)\right)}d{\rm f}_2\eeqa
 {This expression shows} $h_{2, n}(\L_2, \Theta, {\rm G}_2, {\rm g}_2)$ is even in ${\rm g}_2$:
\beqa{parity}
h_{2, n}(\L_2, \Theta, -{\rm g}_2)=h_{2, n}(\L_2, \Theta, {\rm G}_2, {\rm g}_2)\quad \forall\ {\rm g}_2\in \torus\ ,
\eeqa 
so, it affords a Fourier expansion  \beqa{h1nm}h_{2, n}(\L_2, \Theta, {\rm G}_2, {\rm g}_2)=\sum_{m=0}^{+\infty}h_{2, nm}(\L_2, \Theta, {\rm G}_2)\cos{m{\rm g}_2}\eeqa and the claimed expansion~\equ{h1***new}, follows.
 We finally check that $\partial_{{\rm G}_2}h_{1, 0m}\equiv 0$ for all $m\ge 0$. But this is a consequence of the fact that, for   ${\rm r}_1=0$, $h_2$ reduces  to $\frac{1}{2\p}\int_{\torus}\frac{d{\rm l}_2}{\|x_\cK^\ppd\|}$, which is $\G$--independent by Proposition~\ref{prop: kepler}.
Then  the assertion and 
hence the thesis~\equ{Hprop} holds.
We now check that, in the case of $h_2$, one also has $h_{2, nm}=0$ for $n-m$ odd, so, for $n\ge 1$, the stronger identity in~\equ{Hprop1} holds.
Denoting as $c_{np}\in \rational$ the coefficients in the expansion 
$$\cP_n(t)=\sum_{p=0}^n c_{np} t^p$$
where, we recall,  only  $p$'s having the same parity as $n$ appear (an explicit formula for the $c_{np}$'s is available from the first formula in Equation \equ{zeroDer} below), so that
\beqa{Pn1}\cP_n\left(-\sqrt{1-\frac{\Theta^2}{{\rm G}_2^2}}\cos({\rm g}_2+{\rm f}_2)\right)={(-1)^n}\sum_{p=0}^{n}c_{np} \left(1-\frac{\Theta^2}{{\rm G}_2^2}\right)^{p/2}\cos^p({\rm g}_2+{\rm f}_2)\ .\eeqa
 {Using the expansion}
\beqano
\cos^p({\rm g}_2+{\rm f}_2)&=&(\cos{\rm g}_2\cos {\rm f}_2-\sin{\rm g}_2\sin {\rm f}_2)^p\nonumber\\
&=&\sum_{k=0}^p(-1)^k\left(
\begin{array}{llrr}
p\\
k
\end{array}
\right)\sin^k{\rm g}_2\cos^{p-k}{\rm g}_2\sin^k {\rm f}_2\cos^{p-k} {\rm f}_2
\eeqano
and {finally} inserting this expression into~\equ{Pn1} and {afterwards} into~\equ{Qn}, we can write~\equ{Qn} as  a trigonometric polynomial in ${\rm g}_2$ having degree $n$ given by
\beqa{h1n****}
h_{2, n}&=&{(-1)^n}
\sum_{p=0}^{n}
\sum_{k=0}^pc_{np} \hat h_{npk}(\L_2, {\rm G}_2, \Theta)
\sin^k{\rm g}_2\cos^{p-k}{\rm g}_2
\eeqa%
where
\beqano
{\hat h_{npk}(\L_2, {\rm G}_2, \Theta)}&=&{ (-1)^k\left(1-\frac{\Theta^2}{{\rm G}_2^2}\right)^{p/2}
\left(
\begin{array}{llrr}
p\\
k
\end{array}
\right)\frac{1}{2\pi a_2^{n+1}}}\nonumber\\
&&{\int_{\torus}
{\sin^k {\rm f}_2\cos^{p-k} {\rm f}_2}\left(1+\sqrt{1-\frac{{\rm G}_2^2}{\L_2^2}}\cos{\rm f}_2
\right)^{n-1}d{\rm f}_2}
\eeqano
The function under the integral in the expression above has the same parity as $k$, so $\hat h_{npk}$ vanishes for $k$ odd\footnote{{Incidentally, by explicit computation of the integral, we obtain, for even $k$,
\beqano
\hat h_{npk}(\L_2, {\rm G}_2, \Theta)&=& (-1)^k\left(1-\frac{\Theta^2}{{\rm G}_2^2}\right)^{p/2}
\left(
\begin{array}{llrr}
p\\
k
\end{array}
\right)\frac{1}{2\pi a_2^{n+1}}\nonumber\\
&&\sum_{j=0}^{n-1}\sum_{r=0}^{k/2}(-1)^r\left(
\begin{array}{cc}
n-1\\
j
\end{array}
\right)\left(
\begin{array}{cc}
k/2\\
r
\end{array}
\right)\left(1-\frac{{\rm G}_2^2}{\L_2^2}\right)^{j/2}\frac{(j+p-k+2r-1)!!}{(j+p-k+2r)!!}
\eeqano
where only terms with $j$ having the same parity as $p$ appear.
}
}. Therefore, in the summand in~\equ{h1n****} only even indices $k$ appear. But for any even $k$, $\sin^k{\rm g}_2$$\cos^{p-k}{\rm g}_2$ has a Fourier expansion $\sum_{m=0}^p b_m \cos m{\rm g}_2$ where $m$ has the same parity as $p$, which is the same as $n$.  We collect all of the information in the following
 \begin{proposition}\label{HarrK}
 All the assumptions if Theorem~\ref{HP} are verified with $h=h_2$. Therefore, 
the coefficients $ h_{2, nm}$ in  the expansion~\equ{legendre exp}  verify~\equ{Hprop} and, for $n\ge 1$, they verify the stronger identity~\equ{Hprop1}. 
Choosing $\cC=\cK$, the expansion in~\equ{Pn}--\equ{h1nm} holds true, with $h_{2, nm}$ verifying~\equ{Hprop} and~\equ{Hprop1}. In particular, the term $h_{2, 1}$ vanishes identically and $h_{2, 2}$, called {the} dipolar term,  does not depend on ${\rm g}_2$.
 \end{proposition}
 One could ask   {what the last assertion becomes when using}, instead of the $\cK$--map, one of the {more familiar} maps, $\cD$, $\cJ$or $\cP$,  mentioned in the introduction.  As a matter of fact, the same assertion holds, apart {from} parity in the Fourier expansion:

   \begin{proposition}\label{prop: HarrDJP}
Let  ${\rm g}^\cD_2$, ${\rm g}^\cJ_2$ or ${\rm g}^\cP_2$ denote the angles {conjugate} to ${\rm G}_2$,  in the case of the maps $\cD$, $\cJ$ or $\cP$.
In the expansion 
\beqano h_2=\sum_{n=0}^{+\infty}h_{2, n} {\|x^\ppu_\cC\|^n}\qquad \cC=\cD,\ \cJ,\ \cP\eeqano
the coefficients $h_{2, n}$ afford a Fourier expansion $\sum_{m=0}^{+\infty}(a_{nm}\cos{\rm g}^\cC_2+b_{nm}\sin{\rm g}^\cC_2)$, with  $m$ having the parity as $n$ and $a_{nm}$, $b_{nm}$ verifying~\equ{Hprop} and~\equ{Hprop1}. In particular, $h_{1, 1}\equiv0$ and $h_{1, 2}$ does not depend on ${\rm g}^\cD_2$, ${\rm g}^\cJ_2$ or ${\rm g}^\cP_2$, respectively. \end{proposition}
\proof  
The maps  $\cD$, $\cJ$or $\cP$ share the property that $u_1={\rm G}_2$ is one of their impulses. However, the coordinate {conjugate} to ${\rm G}_2$ is different in any of such cases, and is given by the angle, that here we denote as ${\rm g}_{2}^\cD$, ${\rm g}_{2}^g$ or ${\rm g}_{2}^\cP$, formed by a certain ``node'' (we call so a non--vanishing vector in $\real^3$)  with ${\rm P}^\ppd$ in the plane orthogonal to ${\rm C}^\ppd$, with respect to the positive direction determined by ${\rm C}^\ppd$.  The mentioned node is given by:
 $$\n_\cC=\left\{
 \begin{array}{llrr}{\rm k}\times {\rm C}^\ppd\qquad &{\rm if}\ &\cC=\cD\\\\
{\rm C}\times {\rm C}^\ppd &{\rm if}&\cC=\cJ\\\\
{\rm P}^\ppu\times {\rm C}^\ppd &{\rm if}&\cC=\cP
\end{array}
\right.
$$ 
 where ${\rm P}^\ppu$ denotes the direction of the perihelion associated to the Keplerian ellipse of the inner body.
 We then find the following relation
 $${\rm g}_2={\rm g}_{2}^\cC+\f^\cC\qquad \cC=\cD\ ,\ \cJ\ ,\ \cP$$
 where $\f^\cC$ is the angle determined by $\n_\cC$ and $\n_2$ in~\equ{K}. Such function does not depend on ${\rm g}_{2}^\cC$. Since the functions $a_2$, $\Theta$, ${\rm r}_1$, $\zeta_2$, ${\rm f}_2$ in~\equ{Pn}, expressed in terms of $\cD$, $\cJ$, $\cP$, even do not depend on   ${\rm g}_{2}^\cC$, the proof of Proposition~\ref{prop: HarrDJP} follows, replacing such functions into~\equ{Pn}, and using the information given by Proposition~\ref{HarrK}. \qed

\section{Renormalizable integrability}\label{Renormalizable integrability}

{Another  consequence of Theorems~\ref{theo: Harr} and~\ref{partial integral} is that there actually exists a functional dependence between $h_2$ and ${{\rm E}_0}$ {which} we shall   write  explicitly. To this end, we  {premise some abstract consideration}.
 \begin{definition}\label{def: renorm integr}\rm Let $h$, $g$ be two  (commuting) functions
of the form
\beqa{HJ}h(p, q, y, x)=\widehat h({\rm I}(p,q), y, x)\ ,\qquad g(p, q, y, x)=\widehat g({\rm I}(p,q), y, x)\eeqa
where 
\beqa{D}(p, q, y, x)\in \cD:=\cB\times U\eeqa
with $ U\subset \real^2$, $\cB\subset\real^{2n}$ open and connected, $(p,q)=$ $(p_1$, $\cdots$, $p_n$, $q_1$, $\cdots$, $q_n)$   {conjugate} coordinates with respect to the two--form $\o=dy\wedge dx+\sum_{i=1}^{n}dp_i\wedge dq_i$ and ${\rm I}(p,q)=({\rm I}_1(p,q), \cdots, {\rm I}_n(p,q))$, with
$${\rm I}_i:\ \cB\to \real\ ,\qquad i=1,\cdots, n$$
pairwise Poisson commuting:
\beqa{comm}\big\{{\rm I}_i,{\rm I}_j\big\}=0\qquad \forall\ 1\le i<j\le n\qquad i=1,\cdots, n\ .\eeqa
 We say that $h$ is {\it renormalizably integrable via $g$} if there exists a function  $$\widetilde h:\qquad {\rm I}(\cB)\times g(U)\to \real\ , $$ such that
\beqa{renorm}h(p,q,y,x)=\widetilde h({\rm I}(p,q), \widehat g({\rm I}(p,q),y,x))\eeqa
for all $(p, q, y, x)\in \cD$.
\end{definition}}

 \nl
 {\begin{proposition}\label{rem}\rm
If $h$  is renormalizably integrable via $g$, then:\item[{\rm (i)}] 
${\rm I}_1$, $\cdots$, ${\rm I}_n$  are first integrals to $h$ and $g$;
\item[{\rm (ii)}] 
$h$ and $g$ Poisson commute. \end{proposition}
\proof It follows from~\equ{HJ} that
\beqa{comm1}\big\{h, g\big\}=\sum_{1\le i<j\le n}\big\{{\rm I}_i, {\rm I}_j\big\}\big(\partial_{{\rm I}_i}\widehat h\partial_{{\rm I}_j}\widehat g-\partial_{{\rm I}_i}\widehat g\partial_{{\rm I}_j}\widehat h\big)+(\partial_yh\partial_xg-\partial_yg\partial_xh)\eeqa
In this expression,  all the terms  in the summand vanish because of~\equ{comm}, while the last term vanishes because of~\equ{renorm}:
$$\partial_yh\partial_xg-\partial_yg\partial_xh=\partial_g\widetilde h\partial_xg\partial_yg-\partial_yg\partial_g\widetilde h\partial_xg=0\ .$$ This proves (ii).
(i) follows from (ii), replacing the couple $(h, g)$ with $(h, {\rm I}_i)$ or $(g, {\rm I}_i)$, with $i=1$, $\cdots$, $n$. \qed
}

\nl
At level of motion, renormalizable integrability can be rephrased as follows.
\begin{proposition}\label{prop: fixed points}
{L}et $h$ be renormalizably integrable via $g$. Fix a value ${\rm I}_0$ for the integrals ${\rm I}$ and look at the motion of 
 $(y, x)$ under $h$ and $g$, on the manifold ${\rm I}={\rm I}_0$. For any fixed initial datum $(y_0, x_0)$, let 
 $g_0:=g({\rm I}_0, y_0, x_0)$. If
 $\o({\rm I}_0, g_0):=\partial_{g}\tilde h({\rm I}, g)|_{({\rm I}_0, g_0)}\ne 0$, the motion $(y^h(t), x^h(t))$ with initial datum$(y_0, x_0)$ under $h$ is related to the corresponding motion $(y^g(t), x^g(t))$  under  $g$ via
 $$y^h(t)=y^g(\o({\rm I}_0, g_0) t)\ ,\quad x^h(t)=x^g(\o({\rm I}_0,  g_0) t)$$
In particular, under this condition, all the fixed points of $g$
in the plane $(y, x)$ are fixed point to $h$.
Values of $({\rm I}_0, g_0)$ for which $\o({\rm I}_0, g_0)= 0$ provide, in the plane $(y, x)$, curves of fixed points for $h$ (which are not necessarily curves of fixed points to $g$). 
\end{proposition}
\proof All the assertions follow from the formulae, implied  by~\equ{HJ}:
$$\dot y^h=-h_x=-\tilde h_x=-\o({\rm I}_0, g_0) g_x({\rm I}_0, y^h, x^h)$$
and, similarly,
$$\dot x^h=\o({\rm I}_0, g_0) g_y({\rm I}_0, y^h, x^h)\ .$$\qed

\nl
{Below, we prove that, under an additional condition, the converse of Proposition~\ref{rem}  holds true.}
{\begin{theorem}\label{easy lemma}
Let  
$h$, $g$ two  commuting functions of the form~\equ{HJ} on the possibly complex domain $\cD$ as in~\equ{D}, with 
${\rm I}_i$
pairwise Poisson commuting.
For any fixed $c=(c_1,\cdots,c_{n})\in{\rm I}(\cB)$, let $\D_c$ be the set of stationary points  of the function $(y,x)\to g(y,x,c_1,\cdots,c_{n})$, and put $U^*_c:=U\setminus\D_c$ Assume that the set 
 $\cD^*:=\bigcup_{(p,q)\in \cB}\Big\{(p,q)\Big\}\times U^*_{{\rm I}(p,q)}$   has full closure. 
 Then 
 $h$ is renormalizably integrable via $g$. \end{theorem}
  \proof  
We firstly observe that, since $\{h, g\}=\{{\rm I}_i, {\rm I}_j\}=0$ for all $1\le i<j\le n$,  using, as in the proof of Proposition~\ref{rem},   Equation~\equ{comm1},
then
\beqa{comm2}\partial_yh\partial_xg-\partial_yg\partial_xh=\partial_y\widehat h\partial_x\widehat g-\partial_y\widehat g\partial_x\widehat h=0\ .\eeqa
  The assumptions and the Implicit Function Theorem ensure that for any given 
 $c=(c_1,\cdots,c_{n})\in \real^{n}$  in the image of the function $(p, q)\in \cB$ $\to$ $({\rm I}_1$, $\cdots$, ${\rm I}_n)$, 
 and $c_{n+1}$ sufficiently close to in the image of $g(c_1,\cdots,c_{n}, \ovl y,\ovl x)$ where $(\ovl y,\ovl x)\in U^*_c$,
  Equation
 $$g(c_1,\cdots,c_{n}, y,x)=c_{n+1}$$
 can be uniquely solved either with respect to $y$ or $x$, via  suitable functions
 $$y={\rm Y}(c_1,\cdots, c_{n+1}, x)\qquad {\rm or}\qquad x={\rm X}(c_1,\cdots, c_{n+1}, y)\ ,$$
 where ${\rm Y}(c_1,\cdots, c_{n+1}, \cdot)$ is defined on a small neighborhood of $\ovl x$, while ${\rm X}$ $(c_1$, $\cdots$, $c_{n+1}$, $\cdot)$ is defined on a small neighbourhood of $\ovl y$. We now consider
the  function
\beqa{hath1}{\rm h}(c_1\cdots c_{n},c_{n+1}):=\widehat h(c_1,\cdots,c_{n}, {\rm Y}(c_1,\cdots, c_{n+1}, x),x)\eeqa
 and/or the function
\beqa{hath2}{\rm h}'(c_1\cdots c_{n},c_{n+1}):=\widehat h(c_1,\cdots,c_{n},y, {\rm X}(c_1,\cdots, c_{n+1}, y))\ .\eeqa
 We have that ${\rm h}$ is $x$--independent, while ${\rm h}'$ is $y$--independent. Let us check the assertion for ${\rm h}$ (for ${\rm h}'$ is specular). Again by the Implicit Function Theorem:
\beqano
{\rm h}_x&=&\widehat h_y(c_1,\cdots,c_{n}, {\rm Y}(c_1,\cdots, c_{n+1}, x),x){\rm Y}_x(c_1,\cdots, c_{n+1}, x)\nonumber\\
&+&\widehat h_x(c_1,\cdots,c_{n}, {\rm Y}(c_1,\cdots, c_{n+1}, x),x)\nonumber\\
&=&-\widehat h_y(c_1,\cdots,c_{n}, {\rm Y}(c_1,\cdots, c_{n+1}, x),x)
\frac{\widehat g_x(c_1,\cdots,c_{n}, {\rm Y}(c_1,\cdots, c_{n+1}, x),x)}{\widehat g_y(c_1,\cdots,c_{n}, {\rm Y}(c_1,\cdots, c_{n+1}, x),x)}\nonumber\\
&+&\widehat h_x(c_1,\cdots,c_{n}, {\rm Y}(c_1,\cdots, c_{n+1}, x),x)\nonumber\\
&\equiv&0
\eeqano
because of~\equ{comm2}.
  Choosing, for a  fixed $(p,q)\in \cB$, $(y,x)\in U^*_{{\rm I}(p,q)}$,   $c_1={\rm I}_1(p,q)$, $\cdots$, $c_n={\rm I}_n(p,q)$, $c_{n+1}=g(p,q,y,x)$, we have the thesis on the set $\cD^*$. Then, by smooth continuation, the thesis holds on all of $\cD=\cB\times U$.\qed}
  \begin{remark}\label{rem2}\rm
  We observe that the proof is constructive: it provides the function $\widehat h$ via the formulae~\equ{hath1}--\equ{hath2}.
  \end{remark}
  
  \nl
  {In the following, we prove  that $h_2$ is renormalizably integrable via ${\rm E}_0$ as an application of Theorem~\ref{easy lemma}. Afterwards, in Section~\ref{explicitU}, we shall exhibit, explicitly, the relative function $\widetilde h_2$ realizing~\equ{renorm}. In Section~\ref{A curve of fixed points for}, as a counter--example to the last assertion of Proposition~\ref{prop: fixed points}, we shall exhibit a curve of fixed points for $h_{2}$ which is not so for ${\rm E}_0$.}

\paragraph{Application of Theorem~\ref{easy lemma} to $h_2$ and ${\rm E}_0$. } 
We aim to apply Theorem~\ref{easy lemma} to $h_2$ and ${\rm E}_0$.  As in the former section, we use the coordinates $\cK$ defined in~\equ{K}. This map turns to be useful, because the integrals ${\rm I}_1$, $\cdots$, ${\rm I}_5$  are coordinates of such system and hence depend on $(p, q)$ only via one of the $p$'s or one of the $q$'s: see~\equ{nice integrals}. As a first step, we aim to check that $h_2$ and ${\rm E}_0$ have the form in~\equ{HJ}, with   \beqa{integrals}n=3\ ,\quad {\rm I}=({\rm I}_1, {\rm I}_2, {\rm I}_3)=({\rm r}_1, \L_2, \Theta)\ ,\quad y={\rm G}_2\ ,\quad x={\rm g}_2\eeqa  The expression of ${\rm E}_0$ has been given in~\equ{cG***}, so it turns to be as claimed.
The expression of $h_2$ in terms of $\cK$ has been discussed by~\citet{pinzari18a}, and is 

\beqa{belline}
h_2({\rm r}_1,\L_2,\Theta,{\rm G}_2,{\rm g}_2)&=&\frac{1}{2\p}\nonumber\\
&&\int_{\torus}\frac{d{\rm l}_2}{\sqrt{{\rm r}_1^2+2 {\rm r}_1 a_2\varrho_2 \sqrt{1-\frac{\Theta^2}{{\rm G}_2^2}}\cos({\rm g}_2+{\rm f}_2)+a_2^2\varrho_2^2}}
\eeqa
where
$$\varrho_2=\left(1-\sqrt{1-\frac{{\rm G}_2^2}{\L_2^2}}\cos\zeta_2\right)$$
with $\zeta_2$,
 as above,  the eccentric anomaly, and ${\rm f}_2$ the true anomaly, both depending on $(\L_2, {\rm G}_2, {\rm l}_2)$. We observe  that it is possible to have a closed formula for $h_2$, since the integration in $d{\rm l}_2$ can be written explicitly by means of the eccentric anomaly
$$d{\rm l}_2=\varrho_2 d\zeta_2$$
and the true anomaly  ${\rm f}_2$ can be eliminated via the well known relation
$$\varrho_2\cos({\rm g}_2+{\rm f}_2)=\cos{\rm g}_2\left(\cos\zeta_2-\sqrt{1-\frac{{\rm G}_2^2}{\L_2^2}}\right)-\frac{{\rm G}_2}{\L_2}\sin{\rm g}_2\sin\zeta_2\ .$$
Then we rewrite $h_2$ as

 \beqa{h11}
&&h_2({\rm r}_1,\L_2,\Theta,{\rm G}_2,{\rm g}_2)=\frac{1}{2\p}\int_{\torus}\nonumber\\
&&\ \ \frac{\varrho_2 d\zeta_2}{\sqrt{{\rm r}_1^2+2 {\rm r}_1 a_2{
\sqrt{1-\frac{\Theta^2}{{\rm G}_2^2}}
\left(\cos{\rm g}_2\left(\cos\zeta_2-\sqrt{1-\frac{{\rm G}_2^2}{\L_2^2}}\right)-\frac{{\rm G}_2}{\L_2}\sin{\rm g}_2\sin\zeta_2\right)}+a_2^2\varrho_2^2}}\nonumber\\
\eeqa
which is precisely of the form~\equ{HJ}. As a second step, we check that, for any fixed value of the integrals ${\rm I}$ in~\equ{integrals},  the set of fixed points of ${\rm E}_0$ as a function of $({\rm G}_2, {\rm g}_2)$ is at most one--dimensional in the plane $({\rm g}_2, {\rm G}_2)$.
Indeed, equations
$$\arr{\dst \partial_{{\rm G}_2}{{\rm E}_0}=0\\\\
\dst \partial_{{\rm g}_2}{{\rm E}_0}=0}\ ,$$
which read
\beqa{fixed points}\arr{\dst 2{\rm G}_2\left(1-\frac{{\rm m}_2^2{\rm M}_2{\rm r}_1}{2\L_2^2}\frac{\sqrt{1-\frac{\Theta^2}{{\rm G}_2^2}}}{\sqrt{1-\frac{{\rm G}_2^2}{\L_2^2}}}\cos{\rm g}_2+\frac{{\rm m}_2^2{\rm M}_2{\rm r}_1\Theta^2}{2{\rm G}_2^4}\frac{\sqrt{1-\frac{{\rm G}_2^2}{\L_2^2}}}{\sqrt{1-\frac{\Theta^2}{{\rm G}_2^2}}}\cos{\rm g}_2\right)=0\\\\
\dst {\rm r}_1{\sqrt{1-\frac{\Theta^2}{{\rm G}_2^2}}}{\sqrt{1-\frac{{\rm G}_2^2}{\L_2^2}}}\sin{\rm g}_2=0}\eeqa
define an algebraic set in having positive co--dimension. Then Theorem~\ref{easy lemma} applies and we have the following
\begin{proposition}\label{main prop}
$h_2$ is renormalizably integrable via ${\rm E}_0$. Namely,
there exists a function $\widetilde h_2$ such that
$$h_2({\rm r}_1, \L_2, \Theta, {\rm G}_2, {\rm g}_2)=\widetilde h_2\big({\rm r}_1, \L_2, \Theta, {\rm E}_0({\rm r}_1, \L_2, \Theta, {\rm G}_2, {\rm g}_2)\big)\ .$$
\end{proposition}

\nl
In the two following sections we discuss some  insights of dynamical character, related to  the renormalizable integrability of $h_2$.
\subsection{The explicit expression of $\widetilde h_2$}\label{explicitU} 
The function $\widetilde h_2$ in Proposition~\ref{main prop} can be written explicitly, and this is the purpose of this section. Before doing it, let us premise some algebraic consideration.

\begin{definition}[The class $\cH_*$]\label{def: H*}\rm We call {\it class $\cH_*$}   the set of functions of the form
\beqa{H*class}f(a, b, u, v)=\frac{1}{2\p}\int_{\torus}  \frac{P(u c(w)){dw}}{\sqrt{a^2+ 2ab Q(v s(w))+  b^2P(u c(w))^2}}\eeqa
where: $u\to P(u)$, $u\to Q(u)$ are smooth functions for $u=0$;  $P(0)> 0$;  $u\to Q(u)$ is odd;  $c$, $s$ are periodic functions such that there exist two ``symmetries'', i.e., transformations $\s$, $\s':\torus\to \torus$ verifying $|\partial_w\s|=|\partial_w\s'|\equiv 1$ and
$$c\circ\s =c\ ,\quad c\circ\s'=-c\ ,\quad s\circ\s =-s\ ,\quad s\circ\s'=s$$
\end{definition}
 
\nl
Definition 	\ref{def: H*} implies that any $f\in \cH_*$ is  homogeneous of degree $-1$ in $(a, b)$;  even in all of their arguments 
\beqa{parities}
f(-a, b, u, v)&=&f(a, -b, u, v)=f(a, b, -u, v)=f(a, b, u, -v)\nonumber\\
&=&f(a, b, u, v)\qquad \forall\ (a, b, u, v)\eeqa
and, moreover, verifies
\beqa{U***}
f(1, 0, u, v)=f(0, 1, u, v)=1\quad \forall \ (u, v)\ .
\eeqa

\begin{proposition}\label{prop: HR}
All the functions in $\cH_*$ afford a formal series expansion
\beqa{fij} f=\sum_{h, k} f_{hk}(a, b)u^{2h} v^{2k} \eeqa
with
\beqa{fij***}f_{hk}(a, b)=\frac{a^2 b^2p_{hk}(a, b) }{q(a, b)^{\frac{1}{2}+2(h+k)}}\quad {\rm for}\quad (i, j)\in \natural^2\setminus\{(0, 0)\} \eeqa where $q(a, b)$ is a  positive definite quadratic form and $p_{ij}(a, b)$ are polynomials of degree $4(i+j-1)$ with coefficients in ${\mathbb Q}$, even separately in $a$ and $b$. In particular, for any $f\in \cH_*$, there exist
$r$, $s\in{\mathbb Q}$ such that
\beqa{HR} r f_{10}(a, b)+s f_{01}(a, b)\equiv0\quad \forall\ (a, b)\in \real^2\ .\eeqa

\end{proposition}

\begin{remark}\label{comment}\rm 
We call the identity~\equ{HR} {\it generalised Herman resonance} and underline that its validity  is strongly based on the identity~\equ{U***}.  For the averaged Newtonian potential,~\equ{U***} is guaranteed by the Keplerian property (Proposition~\ref{prop: kepler}).
\end{remark}

\proof
Using the formula~\equ{H*class}, it is easy to prove, by induction, that any $f\in \cH_*$
affords an expansion of the kind
\beqano f=\sum_{i, j} \ovl f_{ij}(a, b)u^{i} v^{j} \eeqano
with
\beqa{ovlpij} \ovl f_{ij}(a, b)=\frac{\ovl p_{ij}(a, b) }{q(a, b)^{\su2+i+j}} \eeqa
where $\ovl p_{ij}(a, b)$ are polynomials in $(a, b)$ and
\beqa{qq}q(a, b)=a^2+P(0)^{2} b^2\ .\eeqa
Using the parity of $f$ with respect to all of its arguments, one has, actually, that
$\ovl p_{ij}$'s are even with respect to $a$ and $b$ separately, and vanish if
$i$, $j$ are not both even, so we have an expansion of the form~\equ{fij}, with 
$f_{hk}=\ovl f_{2h, 2k}$.
Furthermore, since  $f$ is homogeneous of degree $-1$, all of its derivatives with respect to $u$ or $v$ are homogeneous of the same degree. Since $q(a, b)$ is homogeneous of degree $2$ (see~\equ{qq}), we have that the $\ovl p_{2h, 2k}$ in~\equ{ovlpij}  are to be homogeneous of degree $4(h+k)$. Finally,  due to~\equ{U***}, $\ovl p_{2h, 2k}(1, 0)=\ovl p_{2h, 2kj}(0, 1)\equiv0$
for all $(h, k)\ne (0, 0)$. Combining this with parity of $\ovl p_{2h, 2k}$ with respect to $a$ and $b$ separately,~\equ{fij***},
$\ovl p_{2h, 2k}(a, b)=a^2b^2 p_{hk}(a, b)$ where $p(h, k)(a, b)$ has degree $4(h+k-1)$. This proves the former assertion. The latter follows from this,  since, when $h+k=1$,
$$f_{10}=\frac{a^2 b^2p_{10}(a, b) }{q(a, b)^{\frac{5}{2}}}\ ,\quad f_{01}=\frac{a^2 b^2p_{01}(a, b)}{q(a, b)^{\frac{5}{2}}} $$
with $p_{10}$ and $p_{01}$ having degree $0$, namely, $p_{10}$ and $p_{01}\in {\mathbb Q}$. So, one can take $r=-p_{01}$, $s=p_{10}$. \qed

\nl
Let us now proceed to write down an  explicit expression of 
function $\widetilde h_2$ in Proposition~\ref{main prop}.
We let

\beqa{U} {\rm U}(a, b,  u ,  v )=\frac{1}{2\p}\int_{\torus}\frac{(1- u \cos w)dw}{\sqrt{a^2+b^2-2b(a v \sin w+b u \cos w)+b^2 u ^2\cos^2w}}\  ; \eeqa
\beqa{EI} {\cal E}(\L_2,{{\rm E}_0})=\frac{\sqrt{\L_2^2-{{\rm E}_0}}}{\L_2}\qquad {\cI}(\L_2,\Theta, {{\rm E}_0})=\frac{\sqrt{{{\rm E}_0}-{\Theta^2}}}{\L_2}\ .
\eeqa

\nl 
Note that
${\rm U}$ is in the class $\cH_*$, with
\beqano
P(u, v)&=&1-u\ ,\quad Q(v)=v\ ,\quad c(w)=\cos w\ ,\quad s(w)=\sin w\ ,\quad \s(w)=-w\nonumber\\
\s'(w)&=&\p-w
\eeqano
We prove that
\begin{proposition}\label{tildeh2**}$\dst \widetilde h_2({\rm r}_1, \L_2, \Theta, {\rm E}_0)={\rm U}({\rm r}_1, a_2, {\cal E}(\L_2,{{\rm E}_0}), {\cI}(\L_2,\Theta, {{\rm E}_0}))$ \end{proposition}

\proof
Reasoning as in the proof of Theorem~\ref{easy lemma} (see Remark~\ref{rem2}), we invert equation 
$${\rm E}_0({\rm r}_1, \L_2, \Theta, {\rm G}_2, {\rm g}_2)=\ovl{\rm E}_0$$
with respect to ${\rm G}_2$ in the complex field, fixing a value of ${\rm g}_2$. We choose ${\rm g}_2=\frac{\p}{2}$, so that $\cos{\rm g}_2=0$ and the the inversion is immediate:
$${\rm G}_2=\sqrt{\ovl{\rm E}_0}\ .$$
Then $\widetilde h_2({\rm r}_1, \L_2, \Theta, \ovl{\rm E}_0)$ is given by
\beqa{tilde h1}\widetilde h_2({\rm r}_1, \L_2, \Theta, \ovl{\rm E}_0)=h_2\left({\rm r}_1, \L_2, \Theta, \sqrt{\ovl{\rm E}_0}, \frac{\p}{2}\right)\eeqa
Using the formula in~\equ{h11}, we obtain

{\small\beqa{tilde h2}&&\widetilde h_2({\rm r}_1, \L_2, \Theta, {\rm E}_0)=\frac{1}{2\p}\int_{\torus}d\zeta_2\nonumber\\
&&\ \frac{1-{\cal E}(\L_2,{{\rm E}_0})\cos\zeta_2}{\sqrt{{\rm r}_1^2+a_2^2-2a_2({\rm r}_1{\cI}(\L_2,\Theta, {{\rm E}_0})\sin\zeta_2+a_2{\cal E}(\L_2,{{\rm E}_0})\cos\zeta_2)+a_2^2{\cal E}(\L_2,{{\rm E}_0})^2\cos^2\zeta_2}}\eeqa}
with $\cE$, $\cI$ as in~\equ{EI}.
\qed

\begin{remark}\label{h2tilde}\rm

Combining  Propositions~\ref{main prop},~\ref{tildeh2**} with~\equ{cG} and the definitions of ${\rm G}_2$ and $\Theta$ in~\equ{K}, we obtain that, for a generic $\cC$ as in~\equ{C},
$$h_2={\rm U}(\|x_\cC^\ppu\|, a_2, {\cal E}_\cC, {\cI}_\cC)\ ,$$
with 
 \beqano
 \cE_\cC&:=&\left(\sqrt{{\rm e}_2^2 +{\rm e}_2\,\frac{x^\ppu\cdot {\rm P}^\ppd}{a_2}}\right)\circ\cC\	 ,	\nonumber\\
\cI_\cC &:=&\left(\sqrt{\frac{\|x^\ppu\|^2\|{\rm C}^\ppd\|^2-(x^\ppu\cdot {\rm C}^\ppd)^2}{\|x^\ppu\|^2\L_2^2} -{\rm e}_2\,\frac{x^\ppu\cdot {\rm P}^\ppd}{a_2}}\right)\circ\cC\ .\eeqano
 In the Section~\ref{mixed average}, we shall use the following consequence of this.
 \begin{proposition}\label{prop: 1 Delau} Let
 \beqano
 \cC_2:\qquad(\L_2,\ell_2,\ovl u,\ovl v)\in \cA\times\torus\times U\to (y^\ppd, x^\ppd)\in (\real^3)^2\eeqano
 where $U$ is a domain of $\real^4$,
 verify~\equ{2B}  and let $\tilde v\in \real^3$. Then
$$\frac{1}{2\p}\int_\torus\frac{d\ell_2}{\|\widetilde v-x^\ppd_{\cC_2}\|}={\rm U}(\| \widetilde v \|, a_2, {\cal E}_2, {\cI}_2)\ ,$$
\beqa{EI1}
&&\cE_2:=\sqrt{{\rm e}_{2,\cC_2}^2 +{\rm e}_{2,\cC_2}\,\frac{\widetilde v\cdot {\rm P}_{\cC_2}^\ppd}{a_2}}\	 ,	\nonumber\\ &&\cI_2:=\sqrt{\frac{\|\widetilde v\|^2\|{\rm C}_{\cC_2}^\ppd\|^2-(\widetilde v\cdot {\rm C}_{\cC_2}^\ppd)^2}{\|\widetilde v\|^2\L_2^2} -{\rm e}_{2, \cC_2}\,\frac{\widetilde v\cdot {\rm P}_{\cC_2}^\ppd}{a_2}}\ ,\eeqa
where the sub--fix $\cC_2$ denotes the composition with $\cC_2$.
\end{proposition}
\proof Choose  $\cC=\id\otimes \cC_2$  in~\equ{C}; namely, 
such  that $u=(\widetilde u, \ovl u)$, $v=(\widetilde v, \ovl v)\in \real^3\times \real^2$, with $(y^\ppu, x^\ppu)\circ\cC=(\widetilde u, \widetilde v)\in \real^3\times \real^2$, and $(y^\ppd, x^\ppd)\circ\cC=(y^\ppd, x^\ppd)\circ\cC_2$, depending only on $(\L_2,\ell_2, \ovl u, \ovl v)$. \qed

\end{remark}

\subsection{A curve of fixed points for $h_2$ (which is not so for ${\rm E}_0$)}\label{A curve of fixed points for}
 
{Proposition \ref{main prop} implies that any level set (in the plane $(\GG, {\rm g})$) to ${\rm E}_0$ is also a level set of $h_2$ and hence, in particular, any fixed point to ${\rm E}_0$ is so to $h_2$.  Here we prove that the converse is not true:
}
\begin{proposition} If $\Theta\ne 0$ and ${\rm r}_1/a_2$ is sufficiently small, in the plane $({\rm G}_2, {\rm g}_2)$, there exists at least a curve  of fixed points of $h_2$ which is a level set of ${\rm E}_0$, but is not a curve a fixed points to it. \end{proposition} \proof In principle, to find any such curve, one should solve equation $\o({\rm I}, \tilde h):=\partial_{g}\tilde h({\rm I}, \tilde h)= 0$. In the case of $h_2$, such equation seems too difficult, so we shall use a perturbative {approach}.
We look  at the Taylor expansion~\equ{Pn} of $h_2$ in~\equ{belline} in powers of  ${\rm r}_1$. Letting $\varepsilon:=\frac{\rm r_1}{a_2}$, we obtain
\beqano
h_2&=&\frac{1}{a_2}\Big[1-\frac{\varepsilon^2}{4}\frac{\L_2^3(3\Theta^2-{\rm G}_2^2)}{{\rm G}_2^5}-\frac{3}{8}\varepsilon^3\sqrt{1-\frac{{\rm G}_2^2}{\L_2^2}}\sqrt{1-\frac{\Theta^2}{{\rm G}_2^2}}\frac{\L_2^5}{{\rm G}_2^5}\Big(1-5\frac{\Theta^2}{{\rm G}_2^2}\Big)\cos{\rm g}_2\nonumber\\
&+&{\rm O}(\varepsilon^4)\Big]\eeqano
By~\equ{tilde h1}, a corresponding expansion for the function $\widetilde h$ in~\equ{tilde h2} is obtained letting ${\rm G}_2=\sqrt{{\rm E}_0}$ and ${\rm g}_2=\frac{\p}{2}$. We obtain:

 $$\widetilde h_2({\rm r}_1,\L_2,\Theta,{{{\rm E}_0}})=\frac{1}{a_2}\Big[1-\frac{\varepsilon^2}{4}\frac{\L_2^3(3\Theta^2-{{\rm E}_0})}{{{\rm E}_0}^{5/2}}+{\rm O}(\varepsilon^4)\Big]$$
We study  Equation
\beqa{linearized equation}\tilde\o=\partial_{{\rm E}_0}\widetilde h=-\frac{\varepsilon^2}{a_2}\left[\frac{\L_2^3}{4}\frac{-15\Theta^2+3{\rm E}_0}{{\rm E}_0^{7/2}}+{\rm O}(\varepsilon^2)\right]=0\eeqa
via the Implicit Function Theorem, for small $\varepsilon$.
Neglecting the ${\rm O}(\varepsilon^2)$ inside parentheses, we obtain the solution
 $${\rm E}_0=5\Theta^2\ .$$
 The non--degeneracy condition at this solution is verified, since indeed
 $$\left.\partial_{{\rm E}_0} \frac{-15\Theta^2+3{\rm E}_0}{{\rm E}_0^{7/2}}\right|_{{\rm E}_0=5\Theta^2}=15\frac{7}{2\Theta^7}-3\frac{5}{2\Theta^7}=\frac{45}{\Theta^2}\ne 0$$
 Then for sufficiently small $\varepsilon$, Equation~\equ{linearized equation} has, {as} a solution, {the following level set of ${\rm E}_0$:}
  $${\rm E}_0=5\Theta^2+{\rm O}(\varepsilon^2)$$
Replacing the formula for ${\rm E}_0$ in~\equ{cG***}, we rewrite such solution as the curve, in the $({\rm G}_2, {\rm g}_2)$ plane,
$$\cS:\quad {\rm G}_2^2-5\Theta^2+{\rm m}_2^2{\rm M}_2{\rm r}_1\sqrt{1-\frac{{\rm G}_2^2}{\L_2^2}}\sqrt{1-\frac{\Theta^2}{{\rm G}_2^2}}\cos{\rm g}_2+{\rm O}({\rm r}_1^2)=0\ .$$
By Proposition~\ref{prop: fixed points}, $\cS$ is a curve  of fixed points for $h_2$. It remains to prove that $\cS$ is not a curve of fixed points   for ${\rm E}_0$. The fixed points of ${\rm E}_0$ are the solutions of the system
\equ{fixed points}. The curve $\cS$ includes a point having coordinates
$${\rm G}_2=\sqrt5\Theta+{\rm O}(\rm r_1^2)\ ,\qquad {\rm g}_2=\frac{\p}{2}+{\rm O}(\rm r_1)$$
which does not solve the system~\equ{fixed points} (it does not solve the second equation). \qed
\begin{remark}
\rm The proof fails for $\Theta=0$, because, in such a case, the leading part in  equation~\equ{linearized equation} has no solution.
\end{remark}

\section{An algebraic property of Legendre polynomials}\label{A algebraic property of Legendre polynomials}

\nl
The Legendre polynomials ${\cal P}_n(t)$, with $\cP_0(t)=1$, $\cP_1(t)=t$, $\cdots$, are defined via the $\varepsilon$--expansion
$$\frac{1}{\sqrt{1-2\varepsilon t+\varepsilon^2}}=\sum_{n=0}^\infty {\cal P}_n(t)\varepsilon^n\ .$$
Many notices on such classical polynomials may be found in~\citet[Appendix B]{giorgilli}.

\nl
The purpose of this section {is to} present an  algebraic property of the $\cP_n$'s.
 Roughly,  it says that a certain average of a Legendre polynomial is still a Legendre polynomial. The author is not aware if it was known before and {if} there is a ``dynamical'' explanation of it.

\begin{lemma}\label{Legendre polynomials1}
Let $t\in \real$, $|t|< 1$, ${\cal P}_n$ the $n^{\rm th}$ Legendre polynomial.
 Then, 
\beq{Q2m}\frac{1}{2\p}\int_{\torus}{\cal P}_n(\sqrt{1-t^2}\cos\theta)d\theta=\d_n{\cal P}_n(t)\eeq
where
$$\d_n=\left\{
\begin{array}{llll}
\dst (-1)^m\frac{(2m-1)!!}{(2m)!!}\qquad&{\rm if}\qquad n=2m\quad{\rm is\ even} \\\\
\dst0&{\rm if}\qquad n\quad{\rm is\ odd}
\end{array}
\right.
$$
\end{lemma}

\nl We shall prove Lemma~\ref{Legendre polynomials1} via the following one.
\begin{lemma}\label{evenLegendre}
The even  Legendre polynomials $P_{2m}(t)$ verify, for any $h=0$, $\cdots$, $m$,
\beqa{zeroDer}
D_\t^hP_{2m}(0)&=&(-1)^{m-h}\frac{h!}{(2h)!}\frac{(2m-2h-1)!!}{(2m-2h)!!}\nonumber\\
D^h_\t P_{2m}(1)&=&
\frac{1}{2^h}\frac{(2m+2h-1)!!}{(2m-1)!!}\frac{(2m)!!}{(2h)!!(2m-2h)!!}\eeqa
where $\t:=t^2$. In particular, the following relation holds
$$(-1)^h\frac{(2h-1)!!}{(2h)!!}D^h_{\t}P_{2m}(0)=(-1)^m\frac{(2m-1)!!}{(2m)!!}D^h_{\t}P_{2m}(1)\ .$$
\end{lemma}
\proof We first  prove the former formula in~\equ{zeroDer}. Let $n\in \natural$, $k=0$, $\cdots$, $n$ with $n-k$ even. We have
$$D_t^k\frac{1}{\sqrt{\varepsilon^2-2t\varepsilon+1}}\Big|_{t=0}=(2k-1)!!\frac{\varepsilon^k}{(1+\varepsilon^2)^{\frac{2k+1}{2}}}\ .$$
Therefore, denoting as $\P_{n}$ the projection over the monomial $\varepsilon^n$,
\beqa{zero der}
D^k_t{\cal P}_n(0)&=&D^k_t\Big(\P_{n}\frac{1}{\sqrt{\varepsilon^2-2t\varepsilon+1}}\Big)\Big|_{t=0}=\P_{n}\Big(D_t^k\frac{1}{\sqrt{\varepsilon^2-2t\varepsilon+1}}\Big|_{t=0}\Big)
=(2k-1)!!\P_{n-k}\frac{1}{(1+\varepsilon^2)^{\frac{2k+1}{2}}}\nonumber\\
&=&\frac{(2k-1)!!}{((n-k)/2)!}D_{\eta}^{(n-k)/2}\frac{1}{(1+\eta)^{\frac{2k+1}{2}}}\Big|_{\eta=0}=(-1)^{(n-k)/2}\frac{(k+n-1)!!}{2^{(n-k)/2}((n-k)/2)!}\nonumber\\
&=&(-1)^{(n-k)/2}\frac{(k+n-1)!!}{(n-k)!!}\ .
\eeqa
Then the desired formula follows, taking $n=2m$, $k=2h$ and noticing  that
$$D^h_\t P_{2m}(0)=\frac{h!}{(2h)!}D^{2h}_tP_{2m}(0)\ .$$

\nl
The proof of the latter  formula in~\equ{zeroDer} is a bit more complicate. We propose an algebraic one.

\nl
 First of all, we change variable
$$t=\sqrt{\t}=\sqrt{1-2z}\ .$$
Since
$$D^h_\t=\frac
{(-1)^h}{2^h}D^h_{z}$$
we are definitely reduced to prove the following identity
\beqa{Dh}
D^h_zP_{2m}(\sqrt{1-2z})\Big|_{z=0}&=&D_{2m, 2h}:=
\frac{(-1)^h}{(2h)!}(2m-2h+2)(2m-2h+4)\cdots (2m)\nonumber\\
&&\times(2m+1)(2m+3)\cdots (2m+2h-1)\ .\eeqa

\nl
To this end, we let
\beqno g(\varepsilon, z):=\frac{1}{\sqrt{\varepsilon^2-2\varepsilon\sqrt{1-2z}+1}}\ ,\eeqno
so that  (analogously to~\equ{zero der}) we may identify
\beqa{projection}
D^h_zP_{2m}(\sqrt{1-2z})\Big|_{z=0}=\P_{2m}D^h_zg(\varepsilon, z)\Big|_{z=0}\ .
\eeqa
We  introduce the auxiliary functions
$$g_{a, b}(\varepsilon, z)=\frac{1}{\big(\varepsilon^2-2\varepsilon\sqrt{1-2z}+1\big)^{\a/2}}\frac{1}{\big(1-2z\big)^{\b/2}}\qquad \a,\ \b\in \real$$
so that $g_{1,0}=g$.
 Observe that the linear space generated by such functions is closed under the derivative operation, since in fact
$$
D_zg_{a, b}(\varepsilon, z)=-\varepsilon\a g_{\a+2,\b+1}(\varepsilon, z)+\b g_{a, b+2}(\varepsilon, z)\ .$$
More in general, by iteration, one finds
\beq{h Der}D^h_zg_{a, b}(\varepsilon, z)=\sum_{j=0}^h c^\pph_{j}\varepsilon^{j}g_{\a+2j, \b+2h-j}(\varepsilon, z)\eeq
where, from the identity
$$D^{h+1}_zg_{a, b}(\varepsilon, z)=D_z\Big(D^{h}_zg_{a, b}\Big)(\varepsilon, z)$$
one easily sees that the coefficients $c^\pph_{j}$, with $j=0$, $\cdots$, $h$ satisfy the  recursion
\beqano
\arr{\dst
c_0^\ppo=1\\\\
\dst c_j^{(h+1)}=-c_{j-1}^{(h)}(\a+2j-2)+(\b+2h-j) c_j^{(h)}\\\\\dst h=0,\ 1,\ \cdots\ ;\qquad  j=0,\ 1,\ \cdots h+1\\\\
\dst c_{-1}^\pph:=0\ ,\qquad c^\pph_{h+1}:=0
}
\eeqano
Let $\ovl c^\pph_j$'s be the numbers defined by
\beqa{ccc}
\arr{\dst
\ovl c_0^\ppo=1\\\\
\dst \ovl c_j^{(h+1)}=-\ovl c_{j-1}^{(h)}(2j-1)+(2h-j) \ovl c_j^{(h)}\\\\\dst h=0,\ 1,\ \cdots\ ;\qquad  j=0,\ 1,\ \cdots h+1\\\\
\dst \ovl c_{-1}^\pph:=0\ ,\qquad \ovl c^\pph_{h+1}:=0
}
\eeqa
corresponding to the case
\beqno\a=1\ ,\qquad \b=0\ .\eeqno
Specialising the formula~\equ{h Der} to this case, we find
\beqano
D^h_zg(\varepsilon, z)\Big|_{z=0}&=&D^h_zg_{1,0}(\varepsilon, z)\Big|_{z=0}\nonumber\\
&=&
\sum_{j=0}^h \ovl c^\pph_{j}\varepsilon^{j}g_{1+2j, 2h-j}(\varepsilon, z)\Big|_{z=0}=\sum_{j=0}^h \ovl c^\pph_{j}\frac{\varepsilon^{j}}{(1-\varepsilon)^{1+2j}}
\eeqano
Therefore, applying~\equ{projection}, we find the desired derivatives
\beqa{derivatives}
D^h_zP_{2m}(\sqrt{1-2z})\Big|_{z=0}&=&\sum_{j=0}^{h}C_{2m, j}\ovl c_j^\pph
\eeqa
with
$$C_{2m, j}:=\frac{(2m-j+1)(2m-j+2)\cdots (2m+j)}{(2j)!}$$

\nl
In order to check~\equ{Dh}, let
$\cP_{2h}(\m)$,  ${\cal Q}_{2j}(\m)$ the   {\sl polynomials} in the {\sl real} variable $\m$ defined as the extensions of $D_{2m, 2h}$, $C_{2m, 2h}$ on the reals, i.e., such that
\beqa{Dh1}
&&\cP_{2h}(2m)=D_{2m, 2h}\ ,\qquad 
{\cal Q}_{2j}(2m)=C_{2m, {2}j}
\eeqa
and let 
$$\cD_{2h}(\m):=\sum_{j=0}^{h}\ovl c_j^\pph{\cal Q}_{2j}(\m)$$
the analogous polynomial extending  the right hand side of~\equ{derivatives}.
We shall prove that
 \beqno \cD_{2h}(\m)=\cP_{2h}({\m})\qquad \forall\ \m\in \real\ ,\quad h=0, 1, \cdots\ ,\eeqno which clearly implies~\equ{Dh}. Note that
$\cD_{2h}(\m)$, $\cP_{2h}(\m)$ have degree $2h$;  $\cP_{2h}(\m)$ vanishes at the odd integers $-(2h-1)$, $-(2h-3)$, $\cdots$, $-1$, and the even integers $0$, $2$, $\cdots$, $2h-2$, while the 
${\cal Q}_{2j}(\m)$'s
have degree $2j$ and vanish at the integers $-j$, $-j+1$, $\cdots$, $j-1$. The last formula in~\equ{Dh1} provides a decomposition of $\cD_{2h}(\m)$ on the basis of the  
${\cal Q}_{2j}$'s. We then do the same for $\cP_{2h}$, i.e., we decompose
\beqano
\cP_{2h}=\sum_{j=0}^{h}
\hat c_j^\pph {\cal Q}_{2j}\ .
\eeqano
We now need to show that 
\beq{deltaj}\hat c^{(h)}_j=\ovl c^\pph_j\qquad \forall\ h=0,\ 1,\ \cdots\ ;\qquad j=0,\ 1,\cdots, h\ .\eeq

\nl
 From the relations 
$$\cP_{2h+2}(\m)=-\frac{(\m-2h)(\m+2h+1)}{2h+2}\cP_{2h}(\m)$$
and 
$$-(\m-2h)(\m+2h+1)=(2h-j)(2h+j+1)-(\m-j)(\m+j+1)$$
the following recursion rule among the coefficients immediately follows
\beqa{cccnew}
\arr{\dst
\hat c_0^\ppo=1\\\\
\dst \hat c^{(h+1)}_j=-\frac{j(2j-1)}{h+1}\hat c^\pph_{j-1}+\frac{4h^2-j^2+2h-j}{2h+2}\hat c^\pph_j\\\\\dst h=0,\ 1,\ \cdots\ ;\qquad  j=0,\ 1,\ \cdots h+1\\\\
\dst \hat c_{-1}^\pph:=0\ ,\qquad \ovl c^\pph_{h+1}:=0
}
\eeqa
Let $$\d^\pph_j:=\hat c^{(h)}_j-\ovl c^\pph_j\ .$$
The formulae in~\equ{ccc} and~\equ{cccnew} imply 
$$\arr{\dst \d^\ppo_0=0\\\\
\dst \d^{(h+1)}_j=-\frac{(2j+1)(j-h-1)}{h+1}\d^\pph_{j-1}+\frac{(2h-j)(j-1)}{2(h+1)}\d^\pph_j\\\\
\dst h=0,\ 1,\ \cdots\ ;\qquad  j=0,\ 1,\ \cdots h+1\\\\
\dst \d_{-1}^\pph:=0\ ,\qquad \d^\pph_{h+1}:=0
}$$
Those relations immediately enforce, by induction, $\d^\pph_j\equiv 0$ for all $h$, $j$, and hence~\equ{deltaj}. \qed

\proof{\it  of Lemma~\ref{Legendre polynomials1}} Let ${\cal Q}_{n}(t)$ denote the left hand side of~\equ{Q2m}. Observe that, {since
any $\cP_n(t)$ has the same parity, in $t$, as $n$ and odd powers of $\cos\theta$ have vanishing average, 
the ${\cal Q}_{2m+1}(t)$'s vanish, while the ${\cal Q}_{2m}(t)$'s are polynomials of degree $m$ in $\t:=t^2$.}
Since also the  even Legendre polynomials $\cP_{2m}$'s are polynomial of degree $m$ in $\t$, we only need to show, e.g., that
$$D^h_{\t}Q_{2m}\big|_{\t=1}=(-1)^m\frac{(2m-1)!!}{(2m)!!}D^h_{\t}P_{2m}\big|_{\t=1}\qquad \forall\ h=0,\ \cdots,\ m\ .$$
The definition of $Q_{2m}$ implies that, for $h=1$, $\cdots$, $m$
$$D^h_{\t}Q_{2m}(1)=(-1)^h\ovl{(\cos\theta)^{2h}}D^h_{\t}P_{2m}(0)=(-1)^h\frac{(2h-1)!!}{(2h)!!}D^h_{\t}P_{2m}(0)\qquad h=0,\ \cdots,\ m$$
where 
$$\ovl{(\cos\theta)^{2h}}:=\frac{1}{2\p}\int_0^{2\p}(\cos\theta)^{2h} d\theta=\frac{(2h-1)!!}{(2h)!!}\ .$$
Using Lemma~\ref{evenLegendre}, we find
$$D^h_{\t}Q_{2m}(1)=(-1)^m\frac{(2m-1)!!}{(2m)!!}D^h_{\t}P_{2m}(1)\qquad h=0,\ \cdots,\ m$$
and hence the thesis follows. \qed

\section{Applications}\label{Applications}
\subsection{An explicit formula for a semiaxes--eccentricities--inclination {expansion} of a ``mixed'' averaged Newtonian potential.  }\label{mixed average}
In this section we assume that the map $\cC$ in~\equ{C} satisfies the following conditions:

\nl
-- the coordinates $(u, v)$ include 
 \beqa{g1}u_1:=\L_1\ ,\quad v_1:=\ell_1\in \torus\ ,\quad v_2:={\rm g}_1\in \torus\eeqa
where, in addition to~\equ{2B}, also the following holds
 \beqa{2B1}
\left(\frac{\|y^\ppu\|^2}{2{\rm m}_1}-\frac{{\rm m}_1{{\rm M}}_1}{\|x^\ppu\|}\right)\circ\cC=
-\frac{{\rm m}_1^3{\rm M}_1^2}{2\L_1^2}=:{\rm h}_{\rm Kep}^\ppu(\L_1)\ ,\eeqa
with suitable other mass parameters ${\rm m}_1$, ${\rm M}_1$; $\ell_1$ in {conjugate} to $\L_1$;

\nl
-- the image of $\cC$ in~\equ{C}
is a domain of  $(y, x)$ where the left hand side of~\equ{2B1} takes negative values;

\nl
--  the instantaneous ellipse ${\mathbb E}_1$ generated by the two--body Hamiltonian~\equ{2B1} has non--vanishing eccentricity;

\nl
-- if ${\rm P}^\ppu$, $\|{\rm P}^\ppu\|=1$ denotes the direction of its perihelion, and, as above,  ${\rm C}^\ppu:=x^\ppu\times y^\ppu$, the angle ${\rm g}_1$ in~\equ{g1} corresponds to the anomaly of ${\rm P}^\ppu$  with respect to a prefixed direction $\n_1$ in (and a prefixed orientation of) the plane orthogonal to ${\rm C}^\ppu$;

\nl
--
$x^\ppd_\cC:=x^\ppd\circ\cC$  and the angle ${\rm f}_1$ (``true anomaly of $x^\ppu_\cC$'')  formed by ${\rm P}_\cC^\ppu$ and $x_\cC^\ppu$ with respect to the orientation established by ${\rm C}^\ppu_\cC$ do not depend on ${\rm g}_1$;

\nl
-- if $\cD_i=(\L_i, {\rm l}_i, p_i, q_i)$, with $p_i=(p_{i1}, p_{i2})\in \real^2$, $q_i=(q_{i1}, q_{i2})\in \real^2$
are the Delaunay coordinates associated to $(y^\ppi, x^\ppi)$,  
and $\cD:= \cD_1\otimes\cD_2:=(\L, {\rm l}, p, q):=(\L_1, \L_2, {\rm l}_1, {\rm l}_2, p_{11}, p_{12}, p_{21}, p_{22}, q_{11}, q_{12}, q_{21}, q_{22})$,
the change of coordinates 
$$\phi_\cC^\cD:\quad \cD_1\otimes\cD_2\to \cC $$
has the form
\beqa{Delchange}\phi_\cC^\cD:\quad\ell_2={\rm l}_2+\f_2(\L, {\rm l}_1, p, q)\ ,\quad (\L, \ell_1, u, v)={\cal F}(\L, {\rm l}_1, p, q)\ .\eeqa

\nl
Our purpose is to provide, under the previous assumptions, a representation formula for the function\footnote{
The reader should not confuse the function $h_{12}$ above with what is commonly called ``doubly averaged Newtonian potential'', defined as
$$\ovl h_{12}:=\frac{1}{(2\p)^2}\int_{\torus^2}\frac{d \ell_1d\ell_2}{\|x_\cC^\ppu-x^\ppd_\cC\|}$$
even though, in the case that ${\mathbb E}_1$ has identically vanishing eccentricity and $\cC$ is regular in this limit, $h_{12}$ and $\ovl h_{12}$ coincide.
} 
$$h_{12}:=\frac{1}{(2\p)^2}\int_{\torus^2}\frac{d{\rm g}_1d\ell_2}{\|x_\cC^\ppu-x^\ppd_\cC\|}$$
which we believe may turn to be useful in applications. We introduce the following
\begin{definition}\label{def: Pi}\rm
For a given power series in the parameter $\varepsilon$
$$g_\varepsilon:=\sum_{n=0}^\infty a_n \varepsilon^n$$
we  denote as $\P_{\varepsilon} g_\varepsilon$  the even power series
$$\P_{\varepsilon} g_\varepsilon:=\sum_{m=0}^\infty (-1)^m\frac{(2m-1)!!}{(2m)!!}{a}_{2m} \varepsilon^{2m}$$
with $(-1)!!:=1$.
\end{definition}
We shall prove the following formula. We let 
 ${\rm U}$ as in~\equ{U} and
\beqa{EI2} {\mathbb E}({\rm r}_1)&=&\sqrt{{\rm e}_{2, \cC}^2 +{\rm r}_1{\rm e}_{2, \cC}\,\frac{{\rm C}_\cC^\ppu\cdot {\rm P}_\cC^\ppd}{\|{\rm C}_\cC^\ppu\|a_{2}}}\nonumber\\
 {\mathbb I}({\rm r}_1)&=&\sqrt{\frac{\|{\rm C}_\cC^\ppu\|^2 \|{\rm C}_\cC^\ppd\|^2-({\rm C}_\cC^\ppu\cdot{\rm C}_\cC^\ppd)^2}{\L_2^2\|{\rm C}_\cC^\ppu\|^2} -{\rm r}_1{\rm e}_{2, \cC}\,\frac{{\rm C}_\cC^\ppu\cdot {\rm P}_\cC^\ppd}{\|{\rm C}_\cC^\ppu\|a_2}}\ ,\eeqa
where the sub--fix $\cC$ denotes the composition with $\cC$.
Then
\begin{proposition}\label{prop: h12}
$h_{12}=\P_{{\rm r}_1}{\rm U}({\rm r}_1, {\rm a}_2, {\mathbb E}({\rm r}_1), {\mathbb I}({\rm r}_1))\big|_{{\rm r}_1=\|x^\ppu_\cC\|}$.\end{proposition}

\begin{remark}[Herman resonance for $h_{12}$]\label{rem: HR}
\rm The functions ${\mathbb E}({\rm r}_1)$, ${\mathbb I}({\rm r}_1)$ in~\equ{EI2} vanish, respectively, in case of zero eccentricity of the exterior planet and mutual inclination. 
Combining Proposition~\ref{prop: HR}, Remark~\ref{comment}, Proposition~\ref{prop: h12}, we obtain
an eccentricity--inclination expansion for $h_{12}$:
\beqano h_{12}=\sum_{h, k} \P_{{\rm r}_1}\left.\left(\frac{{\rm r}_1^2 a_2^2p_{hk}({\rm r}_1, a_2) }{q({\rm r}_1, a_2)^{\frac{1}{2}+2(h+k)}}{\mathbb E}({\rm r}_1)^{2h} {\mathbb I}({\rm r}_1)^{2k}\right)\right|_{{\rm r}_1=\|x^\ppu_\cC\|}\ . \eeqano
The second--order term of this expansion of course exhibits~\equ{HR}, as a byproduct of Proposition~\ref{prop: HR} (because $\P_{\rm r}$ kills the linear  terms in ${\rm r}_1$ in  \equ{EI2} acts on the even terms only modifying the coefficients). This identity reduces to the classical Herman resonance switching to Poincar\'e coordinates with the inner body moving on a circle. In this framework, Herman resonance naturally appears as  a byproduct of parities~\equ{parities}, renormalizable integrability of the Newtonian potential (Proposition~\ref{tildeh2**}), Keplerian property (Proposition~\ref{prop: kepler}) and Lemma~\ref{Legendre polynomials1}. 
\end{remark}


\nl
To prove Proposition~\ref{prop: h12}, we need an equivalent formulation of Lemma~\ref{Legendre polynomials1}, which is as follows.

\begin{proposition}\label{Lemma Legendre}
Let ${\rm r}_1>0$, $\f_1\in \torus$, ${\rm N}^\ppu\in \real^3$, with $\|{\rm N}^\ppu\|=1$, $z^\ppd\in \real^3$,  with $z^\ppd\ne 0$, $z^\ppd\not\parallel {\rm N}^\ppu$. Define $\n:=z^\ppd\times {\rm N}^\ppu$.  Let $z^\ppu({\rm r}_1,\f_1, {\rm N}^\ppu, z^\ppd)$ be such that $z^\ppu\perp {\rm N}^\ppu$, $\|z^\ppu\|={\rm r}_1$ and $\a_{{\rm N}^\ppu}(\n,\ {\rm N}^\ppu\times z^\ppu)=\f_1$. Then, the following  identity holds
\beqa{thesis}
\frac{1}{2\p}\int_{\torus}\frac{d\f_1}{\|z^\ppu({\rm r}_1,\f_1, {\rm N}^\ppu, z^\ppd)-z^\ppd\|}=\frac{1}{{\rm r}_2}\P_{\varepsilon}\frac{1}{\|\varepsilon{\rm N}^\ppu-\widetilde z^\ppd\|}
\Big|_{\varepsilon=\frac{{\rm r}_1}{{\rm r}_2}}
\eeqa
with ${\rm r}_2:=\|z^\ppd\|$,  $\widetilde z^\ppd:=\frac{z^\ppd}{{\rm r}_2}$. Such identity still holds replacing $z^\ppu({\rm r}_1, \f_1, {\rm N}^\ppu, z^\ppd)$ with $z^\ppu({\rm r}_1, \f_1+\hat\f, {\rm N}^\ppu, z^\ppd)$, with any $\hat\f$,  independent of $\f_1$. \end{proposition}

\proof 
Let us decompose $$z^\ppd=(z^\ppd\cdot {\rm N}^\ppu) {\rm N}^\ppu+z^\ppd_\perp$$ 
where
$z^\ppd_\perp:=z^\ppd-(z^\ppd\cdot {\rm N}^\ppu) {\rm N}^\ppu$ is orthogonal to ${\rm N}^\ppu$. Since $z^\ppu$ is orthogonal to ${\rm N}^\ppu$ and $\|z^\ppd_\perp\|=\sqrt{\|z^\ppd\|^2-(z^\ppd\cdot {\rm N}^\ppu)^2}\|={\rm r}_2\sqrt{1-(\widehat z^\ppd\cdot {\rm N}^\ppu)^2}$, we have
$$z^\ppu\cdot z^\ppd=z^\ppu\cdot z^\ppd_\perp=\|z^\ppu\|\| z^\ppd_\perp\|\cos\psi={\rm r}_1{\rm r}_2\sqrt{1-(\widehat z^\ppd\cdot {\rm N}^\ppu)^2}\cos\psi$$
where $\psi$ is the convex angle formed by $z^\ppu$ and $z^\ppd_\perp$. But $\psi$ is related to $\f_1$ via
$$\psi=\|\p-\f_1\|$$
therefore, $\cos\psi=-\cos\f_1$. This readily implies 
\beqa{D**}\|z^\ppu({\rm r}_1,\f_1, {\rm N}^\ppu, z^\ppd)-z^\ppd\|=\sqrt{{\rm r}_1^2+2 {\rm r}_1  {\rm r}_2\sqrt{1-({\rm N}^\ppu\cdot \widehat z^\ppd)^2}\cos\f_1+{\rm r}_2^2}\ .\eeqa
 We now use this in the expansion of the inverse distance
$$\frac{1}{{\rm D}({\rm r}_1, \f_1, {\rm N}^\ppu, z^\ppd)}=\frac{1}{\sqrt{{\rm r}_1^2+2 {\rm r}_1  {\rm r}_2\sqrt{1-({\rm N}^\ppu\cdot \widehat z^\ppd)^2}\cos\f_1+{\rm r}_2^2}}$$
in terms of Legendre polynomials
$$\frac{1}{{\rm D}({\rm r}_1, \f_1, {\rm N}^\ppu, z^\ppd)}=\frac{1}{{\rm r}_2}\sum_{n=0}^\infty (-1)^n\Big(\frac{{\rm r}_1}{{\rm r}_2}\Big)^n {\cal P}_n\Big(\sqrt{1-\frac{(z^\ppd\cdot {\rm N}^\ppu)^2}{{\rm r}_2^2}}\cos\f_1\Big)\ .$$
To conclude, we only need to use Lemma~\ref{Legendre polynomials1}, so  that 
$$\frac{1}{2\p}\int_{\torus}{\cal P}_n\Big(\sqrt{1-\frac{(z^\ppd\cdot {\rm N}^\ppu)^2}{{\rm r}_2^2}}\cos\f_1\Big)d\f_1=\d_n P_{n}\Big(\frac{z^\ppd\cdot {\rm N}^\ppu}{{\rm r}_2}\Big)\ ,$$
which is a rewrite of the thesis. From the formulae from~\equ{D**} on, it follows that the identity~\equ{thesis} still holds replacing $z^\ppu({\rm r}_1, \f_1, {\rm N}^\ppu, z^\ppd)$ with $z^\ppu({\rm r}_1, \f_1+\hat\f, {\rm N}^\ppu, z^\ppd)$, for any $\hat\f$  independent of $\f_1$. \qed

\nl
We can now proceed to prove Proposition~\ref{prop: h12}. We do it in three steps.

\paragraph{1st Step. Application of  {Proposition~\ref{Lemma Legendre}}}  Let $\cC$ and $\n_1\in \real^3\setminus\{0\}$ be as said at the beginning of this section.  As a first step, we aim to compute the ${\rm g}_1$--average applying Proposition~\ref{Lemma Legendre}.  If $${\rm N}^\ppu:=\frac{{\rm C}_\cC^\ppu}{\|{\rm C}_\cC^\ppu\|}\ ,\qquad \n:=x_\cC^\ppd\times {\rm N}_\cC^\ppu$$ then
$$\a_{\rm {\rm N}^\ppu}\big(\n, {\rm N}^\ppu\times x_\cC^\ppu\big)={\rm g}_1+v_1+\frac{\p}{2}-\hat v\quad {\rm where}\quad \hat v=\a_{\rm C^\ppu}\big(
\n_1, \n
\big)\ .$$ 
Hence,   we can write
$$x^\ppu_\cC=z^\ppu\left(\|x^\ppu_\cC\|, \frac{{\rm C}^\ppu_\cC}{\|{\rm C}^\ppu_\cC\|}, {\rm g}_1+v_1+\frac{\p}{2}-\hat v,  x^\ppd_\cC\right)$$
where $z^\ppu$ is as in Proposition~\ref{Lemma Legendre}.
We apply Proposition~\ref{Lemma Legendre} with this $z^\ppu$, $z^\ppd=x^\ppd_\cC$, $\hat\f=v_1+\frac{\p}{2}-\hat v$, which is independent of ${\rm g}_1$, by assumption. We find
$$\frac{1}{2\p}\int_\torus \frac{d{\rm g}_1}{\|x^\ppu_\cC-x^\ppd_\cC\|}=\left.\left.\frac{1}{{\rm r}_2}\P_{\varepsilon}\frac{1}{\left\|\varepsilon \frac{{\rm C}^\ppu_\cC}{\|{\rm C}^\ppu_\cC\|}-\widetilde x_\cC^\ppd\right\|}\right|_{\varepsilon=\frac{\|x^\ppu_\cC\|}{\|x^\ppd_\cC\|}}=\left.\P_{{\rm r}_1}\frac{1}{\left\|{\rm r}_1 \frac{{\rm C}^\ppu_\cC}{\|{\rm C}^\ppu_\cC\|}-x_\cC^\ppd\right\|}
%
\right|_{{\rm r}_1=\|x^\ppu_\cC\|}\right.
$$
with   $\widetilde x_\cC^\ppd:=\frac{{x^\ppd_\cC}}{ \|x^\ppd_\cC\|}$. Now we average with respect to $\ell_2$. We obtain, interchanging $\P_{\rm r_1}$ and $\int_{\torus}d\ell_2$,
\beqa{1st step}\frac{1}{4\p^2}\int_{\torus^2} \frac{d{\rm g}_1d\ell_2}{\|x^\ppu_\cC-x^\ppd_\cC\|}=\frac{1}{2\p}\left.\P_{{\rm r}_1}\int_{\torus}\frac{d\ell_2}{\left\|{\rm r}_1 \frac{{\rm C}^\ppu_\cC}{\|{\rm C}^\ppu_\cC\|}-x_\cC^\ppd\right\|}
%
\right|_{{\rm r}_1=\|x^\ppu_\cC\|}
\eeqa

\paragraph{2nd Step. Switch to Delaunay coordinates.}
We apply $\phi_\cC^\cD$ in~\equ{Delchange} to~\equ{1st step}. We obtain
\beqano
\left(\frac{1}{4\p^2}\int_{\torus^2} \frac{d{\rm g}_1d\ell_2}{\|x^\ppu_\cC-x^\ppd_\cC\|}\right)\circ\phi_\cD^\cC&=&
\frac{1}{2\p}\left.\left(\P_{{\rm r}_1}\int_{\torus} \frac{d\ell_2}{\|{\rm r}_{1, \cC} \frac{{\rm C}^\ppu_\cC}{\|{\rm C}^\ppu_\cC\|}-x^\ppd_\cC\|}\right|_{{\rm r}_1=\|x^\ppu_\cC\|}\right)\circ\phi_\cD^\cC\nonumber\\
&=&\frac{1}{2\p}\left.\P_{{\rm r}_1}\int_{\torus}\frac{d{\rm l}_2}{\left\|{\rm r}_{1, \cD} \frac{{\rm C}^\ppu_\cD}{\|{\rm C}^\ppu_\cD\|}-x_\cD^\ppd\right\|}
%
\right|_{{\rm r}_1=\|x^\ppu_\cD\|}\nonumber\\
&=&{\rm U}(
{\rm r}_{1, \cD}
, a_2, {\cal E}_{2, \cD_2}, {\cI}_{2, \cD_2})
\eeqano
where ${\cal E}_{2, \cD_2}$,  ${\cal I}_{2, \cD_2}$ are as in~\equ{EI1}, with $\cC_2=\cD_2$.
We have that
 ${\rm r}_{1, \cD} \frac{{\rm C}^\ppu_\cD}{\|{\rm C}^\ppu_\cD\|}$ depends only on $\cD_1=(\L_1, {\rm l}_1, p_1, q_1)$, while $x^\ppd_\cD$ depends only on $\cD_2=(\L_2, {\rm l}_2, p_2, q_2)$. We have used Proposition~\ref{prop: 1 Delau} with a given $\widetilde w\in \real^3$, $\cC_2=\cD_2$ and {\it next} we have taken $\widetilde w={\rm r}_{1, \cD} \frac{{\rm C}^\ppu_\cD}{\|{\rm C}^\ppu_\cD\|}$.
\paragraph{3rd Step. }
Applying $(\phi_\cC^\cD)^{-1}$ we conclude the proof. \qed
%

\subsection{Is the two--centre Hamiltonian renormalizably integrable?}

In this section, we outline an underlying open problem {in} the framework of the paper. We pose a conjecture, that we aim to study in further work, which, if proved, may be applied to the two--centre Hamiltonian~\equ{Jc}, so as to obtain a stronger assertion than Proposition~\ref{main prop}.

\nl
Throughout the section, $V\subset \real$, ${\cal U}\subset\real^2$ are domains,  $(I,\f)\in \cI\times \torus$, $(p,q)\in {\cal U}$  are pairwise conjugate canonical coordinates. 
We shall be concerned with real--analytic\footnote{Following the standard terminology, a real function $ h $ is said to be real--analytic on a domain $\cP\subset \real^p$ if there exists an open set $\hat\cP$, with $\cP\subset\hat\cP\subset \complex^p$, such that $ h $ has a holomorphic extension on $\hat\cP$.} functions (``Hamiltonians'') for  $(I, \f, p,q,\m)\in \cP=\cI\times \torus\times {\cal U}\times (-\m_0, \m_0)$ having the form:
\beqa{H} h = h _0(I)+\m f(I,\f, p, q, \m)\ \qquad {\rm with}\qquad  h _0(I)\not\equiv 0\quad {\rm on}\quad V\eeqa
 
 \begin{definition}\rm
We say that $\ovl h $ is  {\it in $p$--normal form} if there exist $\{\ovl h _k(\ovl I, \ovl p,\ovl q)\}_{k=0, \cdots, p}$ $h_k: V\times {\cal U}\to \real$ such that
$$\ovl h(\ovl I, \ovl\f, \ovl p,\ovl q)=\sum_{k=0}^p\ovl h _k(\ovl I, \ovl p,\ovl q)\m^k+{\rm O}(\m^{p+1})\quad \forall (\ovl I, \ovl\f, \ovl p,\ovl q)\in \cP\ .$$

\end{definition}

\nl
The following result is well known and hence will be not discussed.

\begin{proposition}
Let $h$ be as in~\equ{H}.
For any $p\in \natural$ it is possible to find a  real--analytic,  canonical and $\m$--close to the identity transformation
$$\phi:\ (\ovl I, \ovl\f, \ovl p,\ovl q)\in\cP\to (I, \f,  p, q)\in \cP$$
such that $\ovl h := h \circ\phi$ is in $p$--normal form:
$$\ovl h(\ovl I, \ovl\f, \ovl p,\ovl q)=\sum_{k=0}^p\ovl h _k(\ovl I, \ovl p,\ovl q)\m^k+{\rm O}(\m^{p+1})\quad \forall (\ovl I, \ovl\f, \ovl p,\ovl q)\in \cP\ ,$$
with $\ovl h_0(\ovl I, \ovl p,\ovl q)=h_0(\ovl I)$.
 \end{proposition}

 \begin{lemma}\label{trivial lemma}
Let $ h $ be  in $p$--normal form and let $ g $ be a  first integral of $ h $. Then

\nl
{\rm (i)}   $ g $ 
is in $p$--normal form;

\nl
{\rm (ii)} $\{\ovl h _1, \ovl g \}={\rm O}(\m^{p})$, where $\ovl h _1(I,  p, q, \m):=\sum_{k=1}^p\ovl h _k(\ovl I, \ovl p,\ovl q)\m^k$
\end{lemma}
\proof (i) Let
\beqano g (I, \f, p,q,\mu)=\sum_{k=0}^\infty\ovl g _k(I, \f, p, q)\m^k\eeqano
denote  the  Taylor--Maclaurin series in $\m$ of $ g $.
We prove that   the functions $\ovl g _j$ are $\f$--independent for all $0\le j\le p$. We proceed by induction on $j$.
Since $ h (\cdot, \cdot, \cdot, \cdot, \m)$ and $ g (\cdot, \cdot, \cdot, \cdot, \m)$ Poisson--commute for all $\m\in (-\m_0, \m_0)$, we find 
$$\{ h _0, \ovl g _0\}=\partial_I h _0(I)\partial_\f \ovl g _0(I,\f, p, q)\equiv0$$
where we have used  that $ h _0$ depends only on $I$.
Since, by assumption, $\partial_I h _0(I)\not\equiv0$, it follows that $\partial_\f \ovl g _0(I,\f, p, q)\equiv0$ and hence $\ovl g _0(I, \f, p, q)$ is $\f$--independent, and hence $\ovl g _0(I, \f, p, q)=\ovl g _0(I, 0, p, q)= g _0\big(I, (p, q)\big)$ for all $\f\in \torus$, with 
$ g _0\big(I, (p, q)\big)$ as in (ii).
 So the step $j=0$ is proved. Assume now that, for a given $0\le j<p$ and any $0\le k\le j$, $\ovl g _k$ is $\f$--independent. Namely,
 $\ovl g _k(I,\f, p, q)= g _k\big(I, (p, q)\big)$, for some function $ g _k\big(I, (p, q)\big)$, with $0\le k\le j$.
 We prove that $\ovl g _{j+1}$ is so.
 Since $ h $ and $ g $ Poisson commute, \beqa{ovlcomm}\{\ovl h ,\ \ovl g \}={\rm O}(\m^{p+1})\ .\eeqa
 Since $j+1\le p$,  the projection of the left hand side over the monomial $\mu^{j+1}$  vanishes:
\beqano\Big\{ h _0,\ovl g _{j+1}\Big\}+\sum_{k=0}^{j}\Big\{ h _{j-k+1},\  g _{k}\Big\}=0\ .\eeqano
In this identity,  the  term  $\Big\{ h _0,\ovl g _{j+1}\Big\}$ has vanishing $\f$--average, 
because $ h _0$ depends only on $I$, while the term $\sum_{k=0}^{j}\Big\{ h _{j-k+1},\  g _{k}\Big\}$ is $\f$--independent, due to the  fact that the $ h _{j-k+1}$ (by assumption) and the $ g _{k}$ (by the inductive hypothesis) are so. Therefore, such two terms  have to identically vanish separately:
$$\Big\{ h _0,\ovl g _{j+1}\Big\}\equiv0\equiv\sum_{k=0}^{j}\Big\{ h _{j-k+1},\  g _{k}\Big\}\ .$$
The vanishing of the left  hand side implies, as in the base step, that $\ovl g _{j+1}$ is $\f$--independent. The vanishing of the right hand side for all $0\le j+1\le p$ is a rewrite of\footnote{Alternatively, observe that, since $\ovl g $ is independent of $\f$ up to the order ${\rm O}(\m^{p+1})$ and $ h $, $ g $ do Poisson--commute,
by~\equ{ovlcomm},
\beqano\m\{\ovl h _1,\ \ovl g \}=\{ h _0,\ \ovl g \}+\m\{\ovl h _1,\ \ovl g \}+{\rm O}(\m^{p+1})=\{\ovl h ,\ \ovl g \}+{\rm O}(\m^{p+1})=\{ h ,\  g \}+{\rm O}(\m^{p+1})={\rm O}(\m^{p+1})\ .\eeqano} (ii).
\qed

 \begin{corollary}
 For any $p$,
    $$\ovl g^{(p)}_{\rm tr}(I, p, q)=g_0(I, p, q)+\sum_{k=1}^{p-1}\ovl g_k(I, p, q)\m^k\ ,\quad \ovl h ^{(p)}_{1, \rm tr}(I, p, q)=\sum_{k=0}^{p-1}\ovl  h _{k+1}(I, p, q)\m^k$$
    verify
    $$\{\ovl g^{(p)}_{\rm tr},\ \ovl h ^{(p)}_{1, \rm tr}(I, p, q)\}={\rm O}(\m^p)$$
 \end{corollary}
We recall that
$$\ovl h _1(\ovl I, \ovl p,\ovl q, 0)=\frac{1}{2\p}\int_0^{2\p} h _1(\ovl I,\f, \ovl p, \ovl q, 0)d\f$$

\begin{definition}\rm
We shall refer to the formal series $\sum_{k=0}^\infty\ovl h _k(\ovl I, \ovl p,\ovl q)\m^k$ as {\it perturbative series in $\m$ to $ h $}.
\end{definition}

\begin{conjecture}
If  $ h $ {as in \equ{H}} has an independent first integral, its perturbative series converges, as well as the perturbative series to $g$. If $\ovl h$, $\ovl g$ denote the sum of the two series, $\ovl h$ is renormalizably integrable via $\ovl g$. \end{conjecture}


\paragraph{\bf Conflict of interests} The author declares that she has no conflict of interest.



\begin{thebibliography}{}
%
%


\bibitem[Abdullah and Albouy, 2001]{abdullahA01}
K.~Abdullah and A.~Albouy, On a strange resonance noticed by {M}. {H}erman, Regul. Chaotic Dyn., 6(4), 421--432 (2001)

\bibitem[Arnold, 1963]{arnold63}
V.I. Arnold, {S}mall denominators and problems of stability of motion in classical
  and celestial mechanics, Russian Math. Surveys, 18(6), 85--191 (1963)

\bibitem[Bekov and Omarov, 1978]{bekovO78}
A.~A. {Bekov} and T.~B. {Omarov}, Integrable cases of the Hamilton-Jacobi equation and some nonsteady
  problems of celestial mechanics, Soviet Astronomy, 22, 366--370 (1978)


\bibitem[Chierchia and Pinzari, 2011]{chierchiaPi11b}
L.~Chierchia and G.~Pinzari,
The planetary {$N$}-body problem: symplectic foliation, reductions
  and invariant tori, 
Invent. Math., 186(1), 1--77 (2011)


\bibitem[Deprit, 1983]{deprit83}
A.~Deprit,
Elimination of the nodes in problems of {$n$} bodies,
Celestial Mech., 30(2), 181--195 (1983)

\bibitem[F\'ejoz, 2004]{fejoz04}
J.~F{\'e}joz,
D{\'e}monstration du `th{\'e}or{\`e}me d'{A}rnold' sur la
  stabilit{\'e} du syst{\`e}me plan{\'e}taire (d'apr{\`e}s {H}erman),
Ergodic Theory Dynam. Systems, 24(5), 1521--1582 (2004)

\bibitem[F\'ejoz, 2002]{fejoz02}
J.~F{\'e}joz,
Quasi periodic motion in the planar three--body problem,
J. Differential Equations, 183,  303--341 (2002)

\bibitem[F\'ejoz, 2013]{fejoz13}
J.~F{\'e}joz,
On ``{A}rnold's theorem" in celestial mechanics --a summary with an
  appendix on the Poincar{\'e} coordinates,
Discrete and Continuous Dynamical Systems, 3, 3555--3565 (2013)

\bibitem[F\'ejoz and Guardia, 2016)]{fejozG15}
J.~F\'{e}joz and M.~Guardia,
Secular instability in the three-body problem,
Arch. Ration. Mech. Anal., 221(1), 335--362 (2016)

\bibitem[Gallavotti, 1986]{gallavotti86}
G.~Gallavotti, Ph\'enom\`enes critiques, syst\`emes al\'eatoires, th\'eories
  de jauge,
Quasi-integrable mechanical systems, 539--624, Amsterdam, 1986

\bibitem[Giorgilli, 2008]{giorgilli}
A.~Giorgilli, Appunti di {M}eccanica {C}eleste (2008). http://www.mat.unimi.it/users/antonio/meccel/meccel.html

\bibitem[Harrington, 1969]{harrington69}
R.~S. Harrington,
The stellar three-body problem,
Celestial Mech., 1(2), 200--209 (1969)

\bibitem[Jacobi, 1842]{jacobi1842}
C.~G.~J. Jacobi,
 Sur l'{\'e}limination des  {n\oe uds} dans le probl{\`e}me des trois
  corps,
Astronomische Nachrichten, Bd XX, 81--102 (1842)

\bibitem[Laskar and Robutel, 1995]{laskarR95}
J.~Laskar and P.~Robutel,
Stability of the planetary three-body problem. {I}. {E}xpansion of
  the planetary {H}amiltonian,
 Celestial Mech. Dynam. Astronom., 62(3), 193--217 (1995)

\bibitem[Meyer, Palaci\'{a}n and Yanguas, 2018]{palacian18}
K.~R. Meyer, J.~F. Palaci\'{a}n, and P.~Yanguas,
Singular reduction of resonant {H}amiltonians,
Nonlinearity, 31(6), 2854--2894 (2018)

\bibitem[Palaci\'{a}n, Sayas and Yanguas, 2013]{palacianSY12}
J.~F. Palaci\'{a}n, F.~Sayas, and P.~Yanguas,
Regular and singular reductions in the spatial three-body problem,
Qualitative Theory of Dynamical Systems, 12(1), 143--182 (2013)

\bibitem[Palaci\'{a}n, Vanegas and Yanguas, 2017)]{palacianVY17}
J.~F. Palaci\'{a}n, J.~Vanegas, and P.~Yanguas,
Compact normalisations in the elliptic restricted three body problem,
Astrophys Space Sci., 362: 215 (2017)

\bibitem[Pinzari, 2015]{pinzari14}
G.~Pinzari,
Canonical coordinates for the planetary problem,
Acta Applicandae Mathematicae, 137(1), 205--232 (2015)

\bibitem[Pinzari, 2018]{pinzari18a}
G.~Pinzari,
 Exponential stability of Euler integral in the three--body problem,
 arXiv: 1808.07633 (2018)

\bibitem[Pinzari, 2018]{pinzari18}
G.~Pinzari,
Perihelia reduction and global {K}olmogorov tori in the planetary
  problem,
 Mem. Amer. Math. Soc., 255(1218) (2018)


\bibitem[Sundman, 1916]{sundman1916}
K.F. Sundman,
Sur les conditions n{\'e}cessaires et  suffisantes pour la
  convergence du d{\'e}veloppement de la fonction perturbatrice dans le
  mouvement plan,
{\"O}fvers af FVS. f\"orh, 58 A:24 (1916)

\end{thebibliography}


\end{document}